\documentclass[a4paper,11pt]{article}
\usepackage[margin=1.2in]{geometry}
\pdfoutput=1
\usepackage{mathtools}

\usepackage[utf8]{inputenc}
\usepackage{amsfonts,amssymb, bm,amsthm}
\usepackage{graphicx,color,epstopdf}
\usepackage{showlabels} 
\allowdisplaybreaks
\newtheorem{mylemma}{Lemma}[section]
\newtheorem{mytheorem}{Theorem}[section]
\newtheorem{mycorollary}{Corollary}[section]
\newtheorem{myremark}{Remark}

\def\XXint#1#2#3{{\setbox0=\hbox{$#1{#2#3}{\int}$}
    \vcenter{\hbox{$#2#3$}}\kern-.5\wd0}}

\def\rId {{\rm Id}}

\def\dd {''}

\def\hu {\hat{u}}
\def\hf {\hat{f}}

\def\bi{{\bf i}}

\begin{document}
\title{Resonance Frequencies of a Slab with Subwavelength Slits: a 
  Fourier-transformation Approach}
\author{Jiaxin Zhou$^1$ and Wangtao Lu$^2$}
\footnotetext[1]{School of Mathematical Sciences, Zhejiang University, Hangzhou
  310027, China. Email: jiaxinzhou@zju.edu.cn.}
\footnotetext[2]{School of Mathematical Sciences, Zhejiang University, Hangzhou
  310027, China. Email: wangtaolu@zju.edu.cn (corresponding author).
  This author is partially supported by Zhejiang Provisional Funds for Distinguished Young Scholars (LR21A010001).}
\maketitle
\begin{abstract}
  This paper proposes a novel, rigorous and simple Fourier-transformation
  approach to study resonances in a perfectly conducting slab with finite number
  of subwavelength slits of width $h\ll 1$. Since regions outside the slits are
  variable separated, by Fourier transforming the governing equation, we could
  express field in the outer regions in terms of field derivatives on the
  aperture. Next, in each slit where variable separation is still available,
  wave field could be expressed as a Fourier series in terms of a countable
  basis functions with unknown Fourier coefficients. Finally, by matching field
  on the aperture, we establish a linear system of infinite number of equations
  governing the countable Fourier coefficients. By carefully asymptotic analysis
  of each entry of the coefficient matrix, we rigorously show that, by removing
  only a finite number of rows and columns, the resulting principle sub-matrix
  is diagonally dominant so that the infinite dimensional linear system can be
  reduced to a finite dimensional linear system. Resonance frequencies are
  exactly those frequencies making the linear system rank-deficient. This in
  turn provides a simple, asymptotic formula describing resonance frequencies
  with accuracy ${\cal O}(h^3\log h)$. We emphasize that such a formula is more
  accurate than all existing results and is the first accurate result especially
  for slits of number more than two to our best knowledge. Moreover, this
  asymptotic formula rigorously confirms a fact that the imaginary
  part of resonance frequencies is always ${\cal O}(h)$ no matter how we place
  the slits as long as they are spaced by distances independent of width $h$.

\end{abstract}
\section{Introduction}
Electro-magnetic wave scattering problems for optical devices with subwavelength
structures have been extensively studied in recent years
\cite{astlalpal00,chenet13, ebblezwol98, seoetal09, stupodgor10, tak01,
  yansam02}. Distinctive phenomena such as extraordinary optical transmission
and local field enhancement have been experimentally observed: light can be
localized and greatly enhanced near subwavelength apertures or holes. Such
features are vastly demanded in many areas, such as biological sensing and
imaging, microscopy, spectroscopy and communication \cite{sarvig07,lienyllud83}.
The underlying theory of field enhancement, arguably, is largely related to wave
frequency matching some resonance frequency in a scattering problem. Roughly
speaking, a resonance frequency refers to certain complex frequency, at which
the scattering problem allows a nonzero wave field to survive under no external
excitation.

In the past decades, a number of subwavelength structures have been studied to
quantitatively analyze the enhancement of wave field
\cite{babbontri10,bonsta94,bontri10, braholsch20,cladurjoltor06,gaoliyua17,holsch19,joltor06a,joltor06b,joltor08,linshizha20,linzha17,linzha18a,linzha18b,
shivol07,shi10}.
Overall, these works have either numerically illustrated or rigorously proved
the fact that: when wave frequency coincides with real part of some resonance
frequency, the wave field can be enhanced by a factor inversely proportional to
certain power of imaginary part of the resonance frequency. Existing theories
treating resonances can be roughly categorized into two approaches:
boundary-integral approach and matched-asymptotics approach. A representative
work of the first approach is \cite{bontri10} by Bonnetier and Triki. They
studied wave scattering by a perfectly conducting half plane with a
subwavelength cavity and proposed a novel integral-equation technique
incorporated with an operator version of Rouch\'e's theorem to asymptotically
describe resonance frequencies; Green functions of subregions were used to
develop governing integral equations on the aperture of the cavity, which, by
asymptotic analysis of the integral kernels, leads to asymptotic behavior of
resonance frequencies. Following this boundary-integral-equation approach,
Babadjian et al. in \cite{babbontri10} studied resonances by two interacting
subwavelength cavities; later on, Lin and Zhang simplified the analyzing
procedure of \cite{bontri10} and studied resonances by a slab with a single slit
\cite{linzha17} or periodic slits \cite{linzha18a,linzha18b}; Gao et al.
\cite{gaoliyua17} studied resonance frequencies by a rectangular cavity with
different conducting boundaries; recently, Lin et al. \cite{linshizha20} studied
Fano resonances in a slab with a periodic array of two subwavelength slits, and
proved that such a subwavelength structure could support real resonance
frequencies, a.k.a, bound states in the continuum \cite{neuwig29} or embedded
eigenvalues \cite{bonsta94,shivol07,shi10}. Joly and Tordeux
\cite{joltor06a,joltor06b,joltor08} and Clausel et al. \cite{cladurjoltor06}
have used the second approach to study resonances of thin slots. It is
worthwhile to mention a nice work of Holley and Schnitzer \cite{holsch19}, who
used matched asymptotic analysis to get a closed-form of leading term of
resonance frequencies of a slab with a single slit; Brand\~ao et al.
\cite{braholsch20} have recently extended this approach to study resonances in a
slab of finite conductivity with a single slit or a periodic array of slits.

This paper aims to establish a novel, rigorous but much simpler theory to study
resonances in a perfectly conducting slab with a finite number of subwavelength
slits of width $h\ll 1$. Unlike the boundary-integral-equation approach which
relies on subregion Green functions, our theory does not make use of any Green
function, but only relies on Fourier transformations. The underlying motivation
is now that subregion Green functions are basically derived by Fourier
transformations, it is certainly more straightforward to study resonances by
such an approach.

Since regions outside the slits are variable separated, by Fourier transforming
the governing equation, we express field in the outer regions in terms of field
derivatives on the aperture. Next, in each slit where variable separation is
still available, waves could be expressed as Fourier series in terms of a
countable basis functions with unknown Fourier coefficients. Finally, by
matching field on the aperture, we establish a linear system of infinite number
of equations governing the countable unknown coefficients. We note that
\cite{holsch19} has already used a similar approach to establish a linear
system, which, however, presumed certain symmetry of wave field and served
only as a numerical solver when normal incidences of real frequencies are
specified. By carefully asymptotic analysis of each entry of the coefficient
matrix and by retaining entries of the matrix up to only leading algebraic order
of $h$, we rigorously show that, by removing only a finite number of rows and
columns, the resulting principle sub-matrix is diagonally dominant so that the
infinite dimensional linear system can be reduced to a finite dimensional linear
system. Resonance frequencies are exactly those frequencies making the linear
system rank-deficient. This in turn provides an asymptotic formula of resonance
frequencies that is accurate up to the order of ${\cal O}(h^3\log h)$, which is
more accurate than all existing results; to the best of our knowledge, this is
the first accurate result especially for a slab with slits of any finite number
more than two. Moreover, our asymptotic formula rigorously confirms a fact that
imaginary parts of all resonance frequencies are always ${\cal O}(h)$ no matter
how we place the slits as long as they are spaced by distances independent of
$h$. As no Green function is used in our formulation, we expect that such an
approach could be more flexible to study resonances in more complicated and
realistic structures, e.g., a slab of finite/infinite conductivities with single
or periodic slits or with single or periodic holes, which we shall report
elsewhere.

The rest of this paper is organized as follows. In section 2, we present the
Fourier-transformation approach by studying a perfectly conducting slab with a
single slit. In section 3, we extend the approach to study resonances of a slab
with multiple slits. In section 4, we draw the conclusion and present some
potential applications of the current method.
\section{Single slit}
\begin{figure}[!ht]
  \centering
  (a)\includegraphics[width=0.4\textwidth]{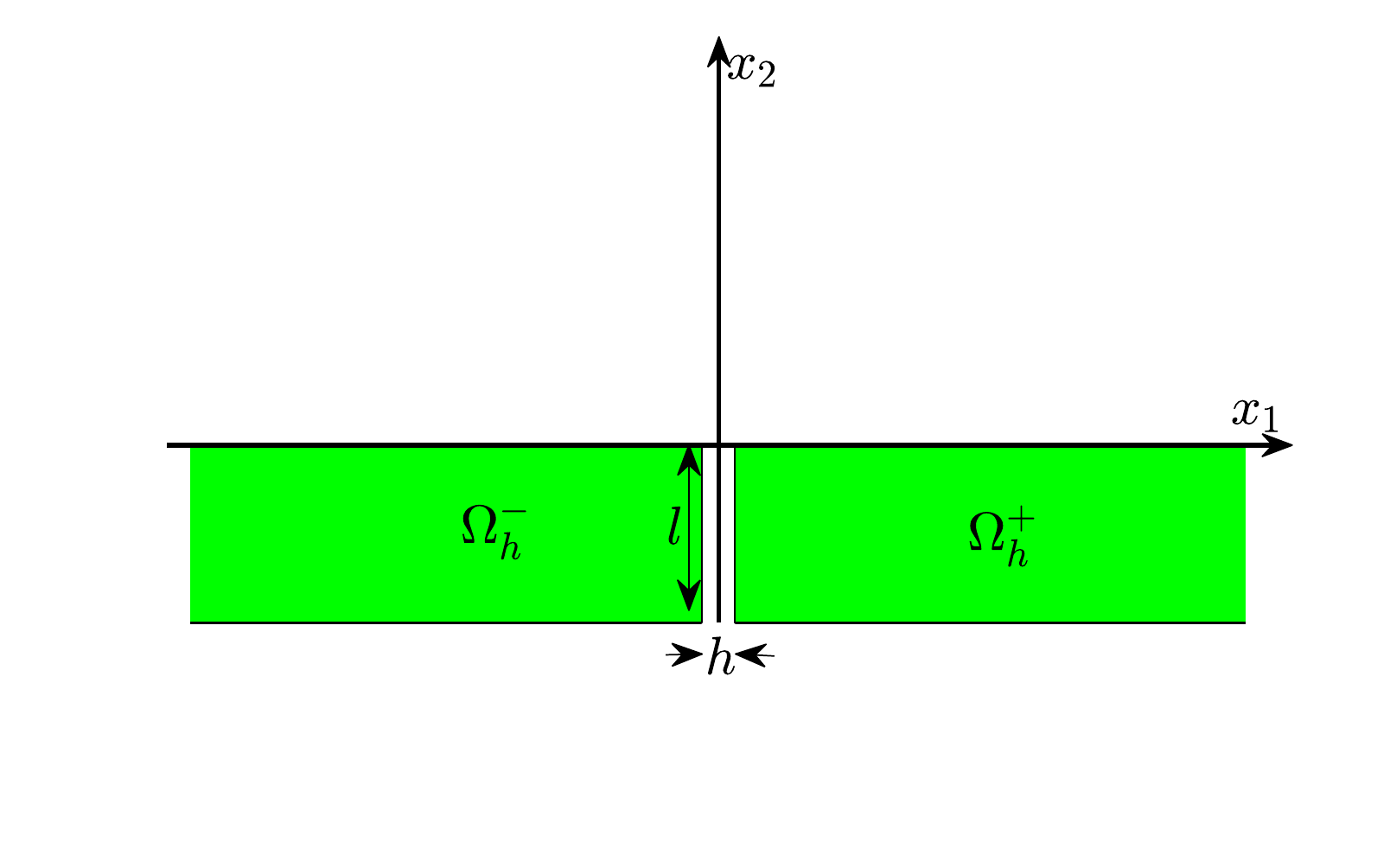}
  (b)\includegraphics[width=0.4\textwidth]{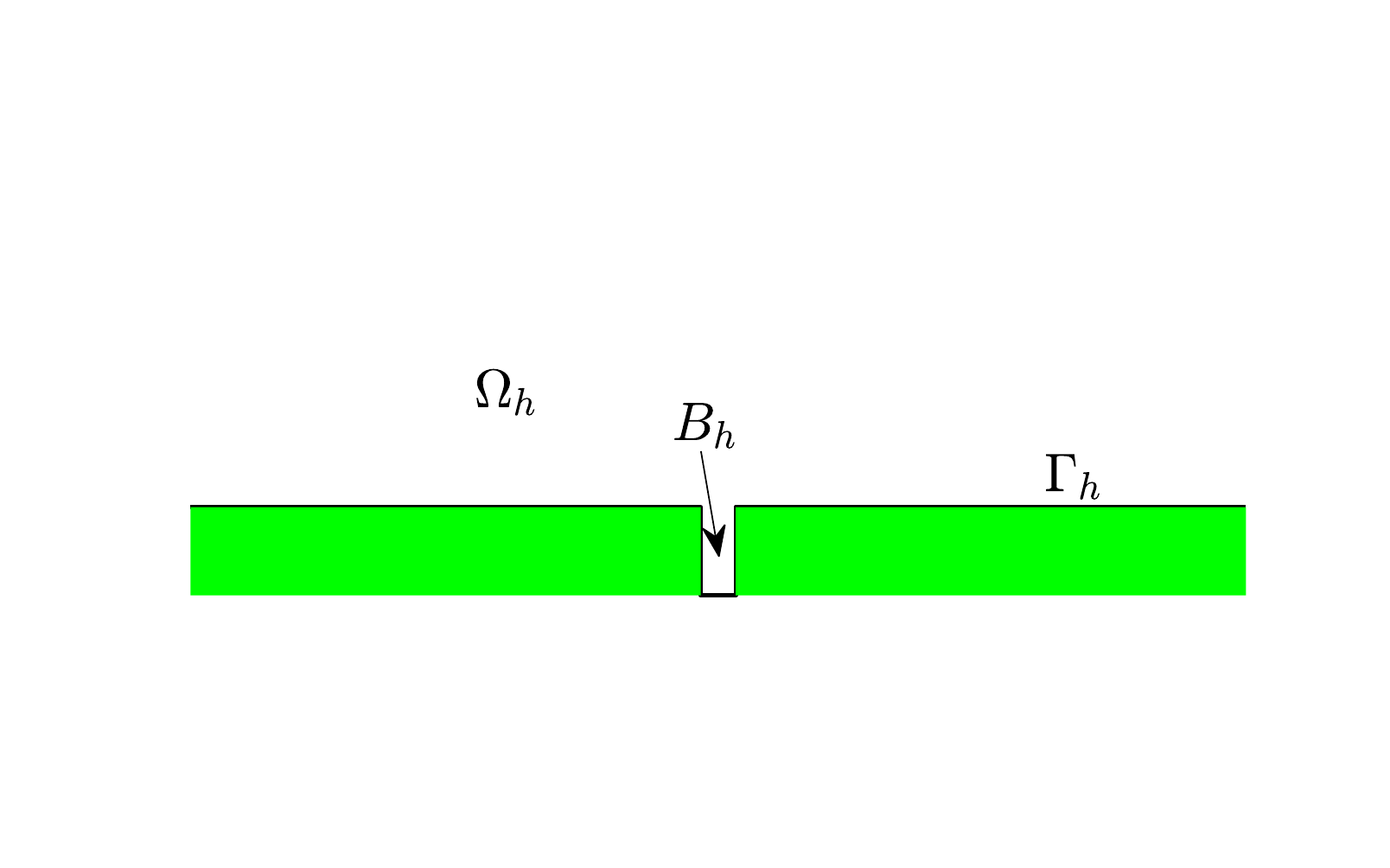}
  \caption{(a). A perfectly conducting slab with a single slit of thickness $l$
    and width $h\ll 1$. (b) The associated half-space problem. }
  \label{fig:model}
\end{figure}

To clarify the basic idea of our Fourier-transformation approach, we begin with
studying resonances of a perfectly conducting slab with a single slit, which has
been studied in \cite{linzha17,holsch19}. Suppose a perfectly conducting slab of thickness $l$ is
perturbed by a slit of width $h\ll 1$, as shown in Fig.~\ref{fig:model} (a).
Then, a TM-polarized
electro-magnetic wave is governed by 
\begin{align}
  \label{eq:helm}
  \Delta u + k^2 u &= 0,\quad{\rm on}\quad \mathbb{R}^2\backslash\overline{\Omega_{h}^{\pm}},\\
  \label{eq:cond}
  \partial_{\nu} u &= 0,\quad{\rm on}\quad \partial \Omega_{h}^+\cup\partial \Omega_{h}^-,
\end{align}
where the two-dimensional Laplace operator
$\Delta=\partial_{x_1}^2+\partial_{x_2}^2$, $k$ is the freespace wavenumber, $u$
denotes the $z$-component of magnetic field, $\Omega_{h}^\pm=\{(x_1,x_2):\pm
x_1>h/2, x_2\in(-l,0)\}$, $\nu$ denotes the outer normal vector along the
boundaries $\partial \Omega_{h}^\pm$. Mathematically, a resonance frequency $k$
refers to a certain value in $\mathbb{C}$, at which there exists a nonzero
solution $u$ solving (\ref{eq:helm}) and (\ref{eq:cond}) and is purely outgoing
at infinity; we shall refer to the nonzero $u$ as a resonance mode in the
following. It is known that when ${\rm Im}(k)>0$, the medium becomes lossy so
that the original problem could not support a nonzero solution. On the other
hand, we expect that ${\rm Im}(k)$ should not be far away from $0$. Thus,
throughout this paper, we shall restrict the searching region to a bounded
domain ${\cal S}=\{k\in\mathbb{C}: {\rm Re}(k)>\epsilon_0>0, |k|<M, \arg(k)\in
(-\frac{\pi}{4},0]\}$ for a sufficiently large constant $M$ and a sufficiently
small constant $\epsilon_0$. By rescaling the variables $x_1$ and $x_2$, we
could assume $l=1$ in the following. Due to the symmetry of the structure, we
split $u$ as the sum of even mode $u^e$ and odd mode $u^o$ about the axis
$x_2=-l/2$ in the following.

\subsection{Even mode}
For simplicity, we suppress the superscript $e$ in this section. Clearly, $u^e$
solves
\begin{align}
  \label{eq:helm:e}
  \Delta u + k^2 u &= 0,\quad{\rm on}\quad \Omega_h,\\
  \label{eq:cond:e}
  \partial_{\nu} u &= 0,\quad{\rm on}\quad \Gamma_h,
\end{align}
where $\Omega_h = \{(x_1,x_2)\in\mathbb{R}^2:x_2>0\}\cup B_h\cup\{(x_1,0):|x_1|<h/2\}$,
$B_h=\{(x_1,x_2):|x_1|<h/2,x_2\in(-l/2,0)\}$, and $\Gamma_h$ denotes the
boundary of $\Omega_h$, as shown in Fig.~\ref{fig:model}(b).

In region $\mathbb{R}_+^2$, Fourier transforming $x_1-$variable in
(\ref{eq:helm:e}) and making use of (\ref{eq:cond:e}), we get
\begin{equation}
  \hu(x_2;\xi) = \frac{\hf(\xi)}{\bi\mu}e^{\bi \mu x_2},
\end{equation}
where $\mu=\sqrt{k^2-\xi^2}$, and 
\begin{align}
  \hu(x_2;\xi) = \int_{-\infty}^{+\infty}u(x)e^{\bi\xi x_1}dx_1,\\
  \label{eq:hf}
  \hf(\xi) = \int_{-h/2}^{h/2}u_{x_2}(x_1,0)e^{\bi\xi x_1}dx_1.
\end{align}
Throughout this paper, unless otherwise specified, we always choose the negative
real axis as the branch cut of $\sqrt{\cdot}$. Thus, by Fourier inverse transform, in $\mathbb{R}_+^2$,
\begin{equation}
  \label{eq:outer}
  u(x) = \frac{1}{2\pi\bi}\int_{-\infty}^{\infty}\frac{\hf(\xi)}{\mu} e^{\bi \mu x_2-\bi\xi x_1}d\xi.
\end{equation}
Inside the slit $B_h$, by method
of variable separations and by (\ref{eq:cond:e}), $u$ can be represented
as
\begin{equation}
  \label{eq:inner}
  u(x) = \sum_{n=0}^{+\infty}b_n\phi_n(x_1)[e^{\bi s_n(x_2+l)} + e^{-\bi s_n x_2}], x_2\in(-l/2,0].
\end{equation}
Here $\{b_n\}$ are unknown scaled Fourier coefficients to be determined,
\begin{align*}
  \phi_n &= \left\{
           \begin{array}{ll}
             \sqrt{\frac{2}{h}}\cos\frac{n\pi x_1}{h} &  n\mid 2;\\
             \bi\sqrt{\frac{2}{h}}\sin\frac{n\pi  x_1}{h} & n\nmid 2,
           \end{array}
           \right.\\
  s_n &= \sqrt{k^2-\left(\frac{\pi n}{h}\right)^2}.
\end{align*}
When $h\ll 1$, using the negative real axis as the branch cut of $\sqrt{\cdot}$
could make $s_n$ a discontinuous and certainly nonholomorphic function for
$k\in{\cal S}$. To resolve this issue, we use the negative imaginary axis as the
branch cut here, only to define $s_n$ \cite{sonlu15}. Clearly, when ${\rm Im}(k)<0$, $s_0=k$
still preserves negative imaginary part while $s_n$ for $n>0$ has positive
imaginary part.

In the above, $\{\phi_n\}$ form a complete and orthogonal basis in the space
$L^2(-h/2,h/2)$ equipped with the natural inner product $(\cdot,\cdot)_{2}$ (they are not orthonormal unless we redefine $\phi_0=\sqrt{1/h}$).
On segment $\Gamma_A:=\{(x_1,0)||x_1|<h/2,x_2=0\}$, we introduce 
the standard Sobolev space $H^{1/2}(\Gamma_A)$ equipped with the following
norm
\[
  ||f||_{H^{1/2}(\Gamma_A)} =\left( \sum_{n=0}^{+\infty}(1+n^2)^{1/2}|\hat{f}_n|^2 \right)^{1/2}, 
\]
where $\hat{f}_n=(f,\phi_n)_2$, and $H^{-1/2}(\Gamma_A)$ as the completion of
$L^2(\Gamma_A)$ w.r.t the following norm: for any $f\in L^2(\Gamma_A)$,
\[
  ||f||_{H^{-1/2}(\Gamma_A)} = \left(  \sum_{n=0}^{+\infty}(1+n^2)^{-1/2}|\hat{f}_n|^2\right)^{1/2}<+\infty.
\]
Furthermore, we define
\[
  \widetilde{H^{-1/2}}(\Gamma_A) = \{f\in H^{-1/2}(\Gamma_A): \exists
  \tilde{f}\in H^{-1/2}(\mathbb{R})\ {\rm such\ that\ } \tilde{f}=0\ {\rm on}\ 
  \mathbb{R}\backslash\overline{\Gamma_A},\ {\rm and}\ \tilde{f}|_{\Gamma_A}=f\}.
\]
Clearly, $\widetilde{H^{-1/2}}(\Gamma_A) = \left( H^{1/2}(\Gamma_A) \right)'$.
Note that the above norms are respectively equivalent to their standard norms
\cite{kre14, mcl00}. Furthermore, let $\ell^2:=\{\{a_n\}_{n=1}^{+\infty}\subset
\mathbb{C}: \sum_{n=1}^{+\infty} |a_n|^2\leq +\infty \}$ be equipped with its
natural norm. We see clearly that $f\in H^{1/2}(\Gamma_A)$ iff
$\{(1+n^2)^{1/4}\hat{f}_n\}\in \ell^2$ and $f\in H^{-1/2}(\Gamma_A)$ iff
$\{(1+n^2)^{-1/4}\hat{f}_n\}\in \ell^2$. In the following, we seek a nonzero
solution $u\in H^{1,{\rm loc}}(\Omega_h)=\{u: u\in H^1(\Omega_h\cap D_R),
D_R=\{x:|x|<R\}\ \forall R>0\}$ solving (\ref{eq:helm:e}) in the
distributional sense. Thus, we require $u|_{\Gamma_A}\in H^{1/2}(\Gamma_A)$ so that the following sequence
\[
  \left\{(1+n^2)^{1/4}b_n(1+e^{\bi s_n l})\right\}\in \ell^2,
\]
and $b_0e^{\bi s_0 l/2}= (u|_{x_2=-l/2},\phi_0)$ with $|b_0|<\infty$. This
implies
\[
  \{a_n:=\sqrt{n}b_{n}\}_{n=1}^{+\infty}\in \ell^2.
\]
On the contrary, if we are given a sequence $\{a_n\}_{n=1}^{\infty}\in \ell^2$
with $b_n=\frac{a_n}{\sqrt{n}}$ and if $|b_0|<\infty$, then equation
(\ref{eq:inner}) defines a solution $u^{\rm in}\in H^1(B_h)$ with
$\partial_{\nu} u^{\rm in}=0$ on
$\partial B_h\cap \Gamma_h$. Thus, for $h\ll 1$,
\begin{align*}
 u_{x_2}^{\rm in}(x_1,0) =  \sum_{n=0}^{+\infty}\bi s_n b_n[e^{\bi s_n l} - 1]\phi_n(x_1)\in \widetilde{H^{-1/2}}(\Gamma_A),
\end{align*}
since 
\begin{align*}
&|(1+n^2)^{-1/4}\bi s_n b_n(e^{\bi s_n l} - 1)|\leq (1+n^2)^{1/4}{\cal O}(b_n(1+e^{\bi s_n l})).
\end{align*}
Such a Neumann data $u_{x_2}^{\rm in}(x_1,0)$ on the real axis defines a unique
solution $u^{\rm ext}\in H^{1}_{\rm loc}(\mathbb{R}_+^2)$ of the Helmholtz
equation (\ref{eq:helm:e}) in $\mathbb{R}_+^2$ with $u^{\rm
  ext}_{x_2}(x_1,0)=u^{\rm in}_{x_2}(x_1,0)$. The two solutions $u^{\rm ext}$
and $u^{\rm in}$ together form a solution $u\in H^{1}_{\rm loc}(\Omega_h)$ of
(\ref{eq:helm:e}) and (\ref{eq:cond:e}) as long as they share the same Dirichlet
data on $\Gamma_A$, i.e.,
\begin{align}
  \label{eq:cont}
  u^{\rm in}(x_1,0) = u^{\rm ext}(x_1,0) = u(x_1,0),\quad |x_1|<h/2.
\end{align}
Thus, by (\ref{eq:hf}) and (\ref{eq:inner}),
\begin{align*}
 \hf(\xi) =&  \sum_{n=0}^{+\infty}\bi s_n b_n[e^{\bi s_n l} - 1]\int_{-h/2}^{h/2}\phi_n(x_1)e^{\bi \xi x_1}dx_1\\
  =&\sqrt{\frac{2}{h}}\sum_{n'=0}^{+\infty}\bi s_{2n'}b_{2n'}[e^{\bi s_{2n'} l} - 1]\frac{2\xi \sin(\xi h/2 + n'\pi)}{\xi^2-\frac{\pi^2 (2n')^2}{h^2}} \nonumber\\
  &+ \sqrt{\frac{2}{h}}\sum_{n'=0}^{+\infty}\bi s_{2n'+1}b_{2n'+1}[e^{\bi s_{2n'+1} l} - 1]\frac{2\xi \sin(\xi h/2 + (n'+1/2)\pi)}{\xi^2-\frac{\pi^2 (2n'+1)^2}{h^2}}\\
  =&\sqrt{\frac{2}{h}}\sum_{n=0}^{+\infty}\bi s_{n}[e^{\bi s_{n} l} - 1]\frac{2\xi \sin(\xi h/2 + n\pi/2)}{\xi^2-\frac{\pi^2 n^2}{h^2}}b_{n}, 
\end{align*}
so that equation (\ref{eq:cont}) implies 
\begin{align*}
  u(x_1,0) =& \sum_{n=0}^{+\infty} \phi_n(x_1) b_n[e^{\bi s_n l} + 1] = \frac{1}{2\pi}\int_{-\infty}^{+\infty}\frac{\hf(\xi)}{\bi\mu}e^{-\bi \xi x_1}d\xi\\
  =&\frac{1}{2\pi}\int_{-\infty}^{+\infty}\frac{1}{\bi\mu}e^{-\bi \xi x_1}\sqrt{\frac{2}{h}}\sum_{m=0}^{+\infty}\bi s_{m}[e^{\bi s_{m} l} - 1]\frac{2\xi \sin(\xi h/2 + m\pi/2)}{\xi^2-\frac{\pi^2 m^2}{h^2}}b_{m}d\xi\\
  =&\sum_{m=0}^{+\infty}b_m[e^{\bi s_{m} l} - 1]\psi_m(x_1),
\end{align*}
where 
\[
  \psi_m(x_1) = \frac{1}{2\pi}\sqrt{\frac{2}{h}}\int_{-\infty}^{+\infty}\frac{s_m}{\mu}\frac{2\xi \sin(\xi h/2 + m\pi/2)}{\xi^2-\frac{\pi^2 m^2}{h^2}}e^{-\bi \xi x_1}d\xi.
\]
Rewriting the above equation in terms of $\{a_n\}\in \ell^2$,
\begin{align}
  \label{eq:gov:an}
  &\phi_0(x_1) b_0[e^{\bi s_0 l} + 1]+\sum_{n=1}^{+\infty} n^{-1/2}\phi_n(x_1) a_n[e^{\bi s_n l} + 1]\nonumber\\
  =&b_0[e^{\bi s_{0} l} - 1]\psi_0(x_1) + \sum_{m=1}^{+\infty}m^{-1/2}a_m[e^{\bi s_{m} l} - 1]\psi_m(x_1).
\end{align}
Consequently, the above arguments in fact imply the following equivalent relation:
\begin{align*}
  &{\rm Finding\ a\ nonzero\ solution}\ u\in H^{1}_{\rm loc}(\Omega_h)\ {\rm of\  
  (\ref{eq:helm:e})\ and\ (\ref{eq:cond:e})}\\
  \Longleftrightarrow &{\rm Finding\ a\ nonzero\ sequence}\ \{a_n\}_{n=1}^{\infty}\in\ell^2\ {\rm and}\ b_0\ {\rm solving\ (\ref{eq:gov:an})}.
\end{align*}

\begin{figure}[!ht]
  \centering
  \includegraphics[width=0.5\textwidth]{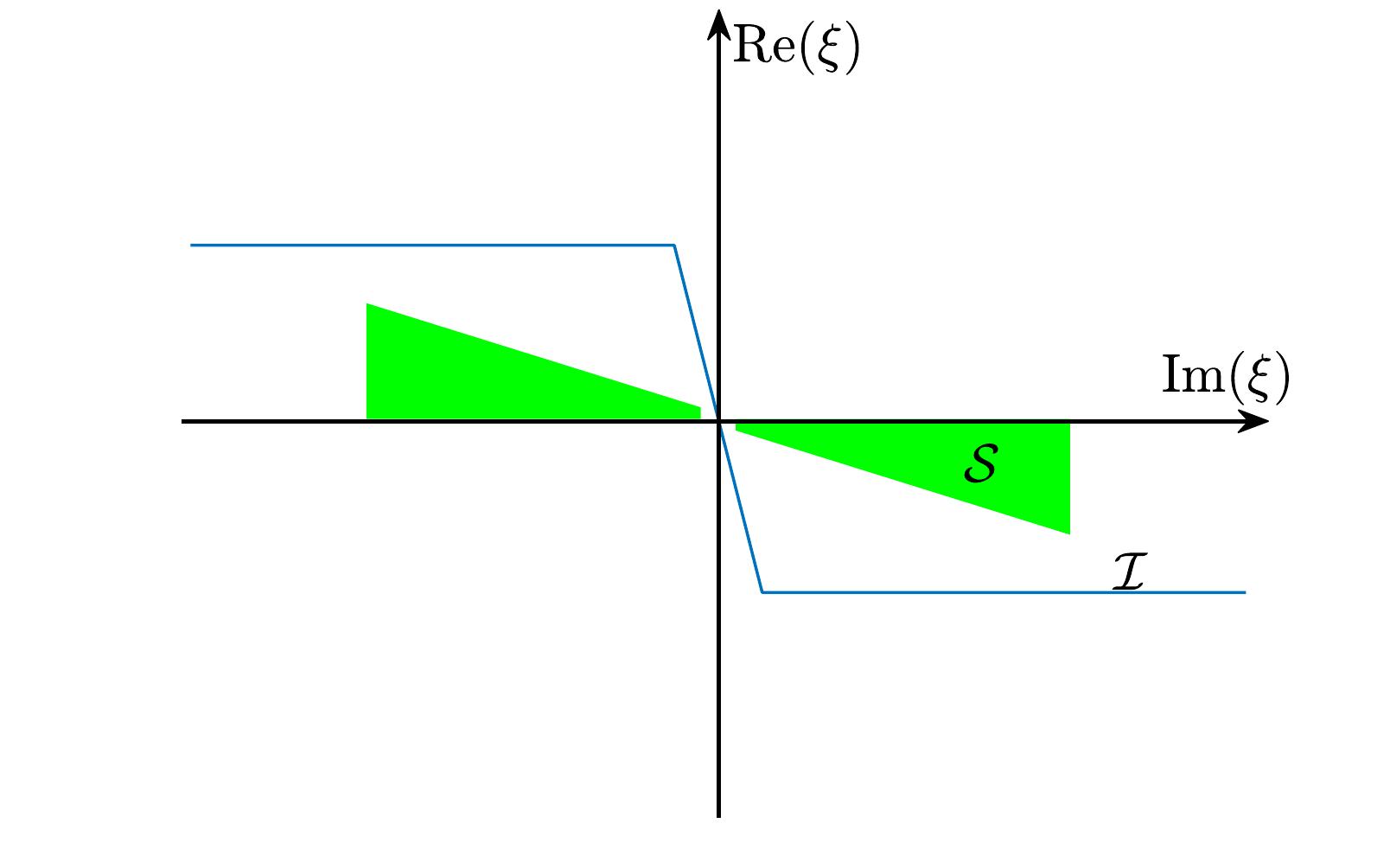}
  \caption{Sommerfeld integral path ${\cal I}$, the searching region $S$ in
    the fourth quadrant and its symmetry. }
  \label{fig:sip}
\end{figure}

Now, taking $\ell^2$-inner product of (\ref{eq:gov:an}) and $\phi_n$ for $n=0,\cdots$, we get
the following linear equations of infinite dimensions:
\begin{align}
  \label{eq:b0}
2b_0[e^{\bi s_0 l} + 1] =& b_0[e^{\bi s_0 l}-1]c_{00} + \sum_{m=1}^{+\infty}a_{2m}c_{2m,0},\\
  \label{eq:bneven}
a_{2n} =& b_0[e^{\bi s_0 l}-1]c_{0,2n}+\sum_{m=1}^{+\infty}a_{2m}c_{2m,2n},\quad n=1,\cdots,\\
  \label{eq:bnodd}
a_{2n-1} =& \sum_{m=1}^{+\infty}a_{2m-1}c_{2m-1,2n-1},\quad n=1,\cdots,
\end{align}
where the sequence $\{d_{mn}\}_{m,n=0}^{\infty}$ is defined as
\begin{align}
  \label{eq:def:dmn}
  d_{mn} =\int_{-\infty}^{+\infty}\frac{1}{\mu h^2}\frac{\xi \sin(\xi h/2 + m\pi/2)}{\xi^2-\frac{\pi^2 m^2}{h^2}}\frac{\xi \sin(\xi h/2 + n\pi/2)}{\xi^2-\frac{\pi^2 n^2}{h^2}} d\xi,
\end{align}
$c_{00} = \frac{4s_0h d_{00}}{\pi}$ and for $m,n\geq 1$,
\begin{align*}
c_{mn} =&\frac{4s_mh}{\pi}\frac{\sqrt{n}(e^{\bi s_m l}-1)}{\sqrt{m}(e^{\bi s_n l}+1)}d_{mn},\\
  c_{m0} =& \frac{4s_mh}{\pi}\frac{(e^{\bi s_m l}-1)}{\sqrt{m}}d_{m0},\\
  c_{0n} =& \frac{4s_0h\sqrt{n}}{\pi}\frac{d_{0n}}{e^{\bi s_n l}+1},
\end{align*}
and we have used the fact that $c_{mn}\equiv 0$ when $m+n\nmid 2$. The above
definition of $d_{mn}$ is well-defined for $k\in\mathbb{R}^+$. However, to
ensure that $d_{mn}$ is a holomorphic function of $k\in{\cal S}$, we should
redefine $d_{mn}$ as follows
\begin{align}
  \label{eq:def:dmn}
  d_{mn} =\int_{\cal I}\frac{1}{\mu h^2}\frac{\xi \sin(\xi h/2 + m\pi/2)}{\xi^2-\frac{\pi^2 m^2}{h^2}}\frac{\xi \sin(\xi h/2 + n\pi/2)}{\xi^2-\frac{\pi^2 n^2}{h^2}} d\xi,
\end{align}
where ${\cal I}$ indicates a Sommerfeld integral path
such that ${\cal S}$ lies above ${\cal I}$ and its symmetry about origin lies
below ${\cal I}$; as shown in Figure~\ref{fig:sip}.

For $h\ll 1$ so that $\epsilon=kh\ll 1$ for $k\in{\cal S}$, we have the following lemma accounting for the asymptotics of $d_{mn}$; its proof is presented in the Appendix.
\begin{mylemma}
  \label{lem:dmn}
  For $\epsilon\ll 1$, the sequence $\{d_{mn}\}_{m,n=0}^{\infty}$ asymptotically behaves as:
  for $k\in\mathbb{Z}$,
  \begin{align}
    \label{eq:dmn}
    d_{mn} = \left\{
    \begin{array}{ll}
      0, & {\rm when}\ m+n\nmid 2;\\
      \frac{\pi}{4} + \frac{\bi}{2}(\gamma-\log 2 - \frac{3}{2}) + \frac{\bi}{2} \log\epsilon + {\cal O}(\epsilon^2\log\epsilon) & {\rm when}\ m=n=0;\\
\bi(-1)^{m/2}C_0(\pi m) + m^{-2}{\cal O}(\epsilon^2\log\epsilon) & 0<m\mid 2, n = 0;\\
\bi(-1)^{n/2}C_0(\pi n) + n^{-2}{\cal O}(\epsilon^2\log\epsilon) & m = 0, 0<n\mid 2;\\
\bi(-1)^{(m-n)/2}C_0^{-}(\pi m,\pi n) -\bi\frac{1}{4m}\delta_{mn} + \frac{\log m -\log n}{m^2-n^2}{\cal O}(\epsilon^2) & 0<m,n\mid 2;\\
\bi(-1)^{(m-n)/2}C_0^{+}(\pi m,\pi n) -\bi\frac{1}{4m}\delta_{mn}+ \frac{\log m -\log n}{m^2-n^2}{\cal O}(\epsilon^2\log\epsilon) & 0<m,n\nmid 2,
    \end{array}
    \right.
  \end{align}
  where $\delta_{mn}$ denotes the Kronecker delta function, for $b,b'\geq \pi$,
\begin{align}
   \label{eq:def:C0}
   C_0(b) &=\int_{0}^{\infty}\frac{1-e^{-kt}}{t(k^2t^2+b^2)}dt\leq \frac{\log b}{b^2},\\
    \label{eq:def:C0-}
    C_0^{-}(b,b')=&\int_{0}^{\infty}\frac{t(1- e^{-t})}{(t^2+b^2)(t^2+b'^2)}dt\leq 
                    \frac{\log b - \log b'}{b^2 - b'^2},\\
    \label{eq:def:C0+}
    C_0^{+}(b,b')=&\int_{0}^{\infty}\frac{t(1+ e^{-t})}{(t^2+b^2)(t^2+b'^2)}dt\leq 
                    \frac{\log b - \log b'}{b^2 - b'^2} + \frac{1}{b^2b'^2},
  \end{align}
  limit is considered when $b=b'$, and the invisible constants in ${\cal
    O}$ notations are independent of $m$ and $n$.
\end{mylemma}
Moreover, as $h\to 0$,
\begin{align*}
  s_mh &= \pi m\bi + {\cal O}(\epsilon^2 m^{-1}),\\
  1\pm e^{\bi s_ml} &= 1 + {\cal O}(e^{-m\pi l/h}).
\end{align*}
 We get the asymptotic expansions  of $c_{mn}$ as follows.
\begin{mylemma}
  \label{lem:cmn}
  For $\epsilon\ll 1$, the sequence $\{c_{mn}\}_{m,n=0}^{\infty}$ asymptotically behaves as:
  for $k\in\mathbb{Z}$,
  \begin{align}
    \label{eq:cmn}
    c_{mn} = \left\{
    \begin{array}{ll}
      0, & m+n\nmid 2;\\
      \epsilon + \frac{2\bi}{\pi}(\gamma-\log 2 - \frac{3}{2})\epsilon + \frac{2\bi}{\pi}\epsilon \log\epsilon + {\cal O}(\epsilon^3\log\epsilon) & m=n=0;\\
4\sqrt{m}(-1)^{m/2}C_0(\pi m) + m^{-3/2}{\cal O}(\epsilon^2\log\epsilon) & n=0<m\mid 2;\\
\frac{4\sqrt{n}}{\pi}\epsilon\bi(-1)^{n/2}C_0(\pi n) + n^{-3/2}{\cal O}(\epsilon^3\log\epsilon) & m = 0<n\mid 2;\\
p_{mn}^{(e)} -\delta_{mn} + \frac{\sqrt{mn}(\log m -\log n)}{m^2-n^2}{\cal O}(\epsilon^2) & 0<m,n\mid 2;\\
p_{mn}^{(o)} -\delta_{mn} + \frac{\sqrt{mn}(\log m -\log n)}{m^2-n^2}{\cal O}(\epsilon^2\log\epsilon) & 0<m,n\nmid 2,
    \end{array}
    \right.
  \end{align}
  where we have defined
  \begin{align}
    \label{eq:pmne}
    p_{mn}^{(e)}=&4\sqrt{mn}(-1)^{(m-n)/2}C_0^{-}(\pi m,\pi n)\\ 
    \label{eq:pmno}
    p_{mn}^{(o)}=&4\sqrt{mn}(-1)^{(m-n)/2}C_0^{+}(\pi m,\pi n). 
  \end{align}
\end{mylemma}

Now for any integer $n\geq 1$, let $A_n^{(e)}$ and $A_n^{(o)}$ be defined as: for any
$\{f_j\}_{j=1}^{+\infty}\in \ell^2$,
\begin{align}
  A_n^{(e)}\{f_j\}=&\{\chi_n(i)\sum_{j=1}^{\infty}(c_{2i,2j} + \delta_{ij})f_{j}\chi_{n}(j)\}_{i=1}^{\infty},\\
  A_n^{(o)}\{f_j\}=&\{\chi_n(i)\sum_{j=1}^{\infty}(c_{2i-1,2j-1}+ \delta_{ij})f_{j}\chi_{n}(j)\}_{i=1}^{\infty},
\end{align}
where 
\[
  \chi_{n}(i) = \left\{
    \begin{array}{ll}
      1 & i\leq n;\\
      0 & {\rm otherwise}.
    \end{array}
  \right.
\]
We have the following theorem.
\begin{mytheorem}
  For $\epsilon\ll 1$, the operators $\{A_n^{(l)}\}_{n=1}^{\infty}, l=e,o$ mapping from
  $\ell^2$ to $\ell^2$ are uniformly bounded, i.e.,
\[
  ||A_n^{(l)}||\leq \frac{1}{2} + \frac{4}{\pi^4} + {\cal O}(\epsilon^2).
\]
As a consequence, there exists a bounded and contracting  operator ${\cal A}^{(l)}: \ell^2\to \ell^2$ such
that for any $\{f_j\}_{j=1}^{\infty}\in \ell^2$,
\[
  {\cal A}^{(e)}\{f_j\} = \{\sum_{j=1}^{+\infty}c_{2i,2j}f_{j}\}_{i=1}^{+\infty}\in \ell^2,\quad
  {\cal A}^{(o)}\{f_j\} = \{\sum_{j=1}^{+\infty}c_{2i-1,2j-1}f_{j}\}_{i=1}^{+\infty}\in \ell^2,
\]
and
\begin{equation}
  \label{eq:Ai:norm}
  ||{\cal A}^{(l)}||\leq \frac{1}{2}+\frac{4}{\pi^4}+{\cal O}(\epsilon^2),
\end{equation}
for $l=e,o$.
\begin{proof}
  Here, we prove the property of ${\cal A}^{(e)}$ only, and shall suppress the
  superscript for simplicity. We first prove the contraction of $A_n$.
  For $\epsilon\ll 1$, we have
  \[
    c_{2i,2j} + \delta_{ij} =8\bi\sqrt{ij}d_{2i,2j} + i^{-3/2}j^{1/2}\frac{\log i-\log
      j}{i^2-j^2}{\cal O}(\epsilon^2),
  \]
  where limit is considered when $i=j$. Note that we should not use the expansion
  (\ref{eq:cmn}) since the neglected part is not symmetric. Let
  $P_n$ and $Q_n$ be operators defined as $A_n$ but with $a_{2i,2j}$ replaced
  respectively by
  \[
    P_{ij}=8\bi\sqrt{ij}d_{2i,2j}:=P_{ij}^{(1)}+\bi P_{ij}^{(2)},\quad
    Q_{ij}=\epsilon^{-2}(c_{2i,2j} + \delta_{ij} - 8\bi\sqrt{ij}d_{2i,2j}),
  \]
  for $i,j\geq 1$. Then 
  \begin{align*}
    ||(Q_{ij})_{n\times n}||_{\rm FRO}^2 \lesssim & \sum_{i,j=1}^{n}i^{-3}j\frac{(\log i-\log j)^2}{(i^2-j^2)^2}<\infty,
  \end{align*}
  where $||\cdot||_{\rm FRO}$ is the Frobenius norm. Since $||(Q_{ij})_{n\times
    n}||_{\rm FRO}$ is strictly increasing w.r.t $n$, and the 2-norm
  $||\cdot||_2\leq ||\cdot ||_{\rm FRO}$, we see clearly that
  $\{Q_n\}_{n=1}^{+\infty}$ is a Cauchy sequence in ${\cal L}(\ell^2;\ell^2)$,
  converging to a bounded operator ${\cal Q}: \ell^2\to \ell^2$. As for the $n\times
  n$ symmetric matrix $(P_{ij}^{(i)})_{n\times n},i=1,2$, its $2$-norm is
  exactly the magnitude of its largest eigenvalue. Thus, we choose to estimate
  the eigenvalue of matrix $(\sqrt{i}P_{ij}^{(i)}\sqrt{j}^{-1})_{n\times n}$,
  which is similar to $(P_{ij}^{(i)})_{n\times n}$. By Lemma~\ref{lem:cmn},
  \begin{align*}
    P_{ij}^{(1)} =& p_{2i,2j}^{(e)} + \frac{\sqrt{ij}(\log i-\log j)}{i^2-j^2} {\cal O}(\epsilon^2) = \frac{\sqrt{ij}(\log i-\log j)}{i^2-j^2}(1+{\cal O}(\epsilon^2)),\\
    P_{ij}^{(2)}=& \frac{\sqrt{ij}(\log i-\log j)}{i^2-j^2}{\cal O}(\epsilon^2),
  \end{align*}
  and 
  \begin{align*}
    \sum_{j=1}^{n}\left|\sqrt{i}\frac{\sqrt{ij}(\log i-\log j)}{i^2-j^2}\sqrt{j}^{-1}\right|\leq& 8\sum_{j=1}^{n}i\frac{\log i-\log j}{4\pi^2(i^2-j^2)}\\
    \leq &\frac{2}{\pi^2}\int_{0}^{+\infty}i\frac{\log x - \log i}{x^2-i^2}dx =\frac{2}{\pi^2}\int_{0}^{+\infty}\frac{\log x}{x^2-1}dx\\
    \leq& \frac{1}{2},
  \end{align*}
  where we notice that the terms in the summation is decreasing in $j$.
  Therefore, for all $n\in\mathbb{N}$, 
  \[
    ||{\rm Re}(P_n)||=||(P_{ij}^{(1)})_{n\times n}||_2\leq
    ||(i^{1/2}P^{(1)}_{ij}j^{-1/2})_{n\times n}||_{1}\leq \frac{1}{2}(1+{\cal O}(\epsilon^2)),
  \]
  for $\epsilon\ll 1$. One similarly obtains that $||{\rm Im}(P_n)||={\cal
    O}(\epsilon^2)$ so that $||P_n||\leq 1/2+{\cal O}(\epsilon^2)$ for all $n$,
  where the invisible constant in the ${\cal O}$-notation is independent of $n$.
  Suppose $\ell^2_{\rm comp}=\{\{f_i\}\in \ell^2: \exists N>0, f_i\equiv 0\ {\rm for}\
  i\geq N\}$. Clearly, $\ell^2_{\rm comp}$ is dense in $\ell^2$. Now, we define ${\cal P}^{(e)}:l_{\rm comp}^2\to \ell^2=(\ell^2)'$ as follows: for any
  $\{f_i\}_{i=1}^{\infty}, \{g_i\}_{i=1}^{\infty}\in \ell^2_{\rm comp}$,
  \begin{equation}
    \label{eq:Pe}
    <{\cal P}^{(e)}\{f_i\},\{g_i\}>:=\sum_{i,j=1}^{\infty}g_jp_{2i,2j}^{(e)}f_i.
  \end{equation}
  Clearly, the above summation is finite since $g_i=f_i\equiv 0$ for
  $i\geq N$. Then, by
  \[
    |<{\cal P}^{(e)}\{f_i\},\{g_i\}>|\leq ||(p_{2i,2j}^{(e)})_{N\times N}||_{2}||\{f_i\}||_{\ell^2}||\{g_i\}||_{\ell^2}\leq \frac{1}{2}||\{f_i\}||_{\ell^2}||\{g_i\}||_{\ell^2},
  \]
  we get from continuous extension theorem that ${\cal P}^{(e)}\{f_i\}\in (\ell^2)'$, and 
  \[
    ||{\cal P}^{(e)}\{f_i\}||_{ \ell^2 }=||{\cal P}_D\{f_i\}||_{( \ell^2 )'} \leq \frac{1}{2}||\{f_i\}||_{\ell^2},
  \]
  which states that ${\cal P}^{(e)}$ is bounded from $\ell^2_{\rm comp}\to \ell^2$. By
  continuous extension theorem again, we see that ${\cal P}^{(e)}$ can be uniquely
  extended as a bounded operator from $\ell^2$ to $\ell^2$ such that
  \[
    ||{\cal P}^{(e)}||\leq \frac{1}{2}.
  \]
  One similarly proves the existence of ${\cal P}_\epsilon: \ell^2\to \ell^2$ defined
  by the elements $\{P_{ij} - p_{2i,2j}^{(e)}\}$ with $||{\cal P}_{\epsilon}||\leq {\cal O}(\epsilon^2)$. and the proof is completed by
  observing that ${\cal A} = {\cal P}^{(e)} + {\cal P}_{\epsilon} + {\cal
    Q}\epsilon^2$.
\end{proof}
\end{mytheorem}
We are ready to solve the infinite dimensional linear system
(\ref{eq:b0}-\ref{eq:bnodd}), which can be restated as: Seek nonzero
$\{a_{j}\}_{j=1}^{\infty}\in \ell^2$ and $|b_0|<\infty$, such that
\begin{align}
  \label{eq:op1}
  2b_0(e^{\bi s_0 l}+1) =& b_0(e^{\bi s_0 l} - 1)c_{00} + <\{a_{2m}\},\{c_{2m,0}\}>_{\ell^2},\\
  \label{eq:op2}
  2\{a_{2n}\} =& b_0(e^{\bi s_0 l} - 1)\{c_{0,2n}\} + {\cal A}^{(e)}\{a_{2n}\},\\
  \label{eq:op3}
  2\{a_{2n-1}\} =& {\cal A}^{(o)}\{a_{2n-1}\},
\end{align}
where we have used the fact that $\{c_{2m,0}\}, \{c_{0,2m}\}\in \ell^2$. As ${\cal
  A}^{(i)}, i=e,o,$ is contracting for $\epsilon\ll 1$, $2{\rm Id}-
{\cal A}^{(i)}$ is invertible, where Id stands for the identity operator.
Consequently, we arrive at our first theorem.
\begin{mytheorem}
  \label{thm:evenres}
 For $h \ll 1$, the system (\ref{eq:op1}-\ref{eq:op3}) has a nonzero solution if and only if
 $k$ solves
 \begin{align}
   \label{eq:gov:k0:e}
   2(e^{\bi k l} + 1) = (e^{\bi k l}-1)\left[  c_{00} + <(2{\rm Id}-{\cal A}^{(e)})^{-1}\{c_{0,2m}\},\{c_{2m,0}\}>_{\ell^2}\right].
 \end{align}
 In fact, the solutions (resonance frequencies) to (\ref{eq:gov:k0:e}) are
 \begin{align}
   \label{eq:kl:e}
   kl=&k_{m,e}-\bi \left[ 1 + \frac{2}{\pi} h  \right]\Delta(\epsilon_{m,e}) - \left[k_{m,e}^{-1}-5\pi^{-1}k_{m,e}^{-1}h\right]\Delta^2(\epsilon_{m,e}) + \bi\left[  k_{m,e}^{-2}-\frac{1}{12}\right]\Delta^3(\epsilon_{m,e}) \nonumber\\
  &+ {\cal O}(\epsilon_{m,e}^3\log\epsilon_{m,e}),\quad m=1,2,\cdots,
 \end{align}
 where $k_{m,e}=(2m-1)\pi$ are the Fabry-P\'erot frequencies and
 $\epsilon_{m,e} =k_{m,e} h\ll 1$,
 \begin{align}
   \label{eq:f}
   \Delta(\epsilon)=&\epsilon + \frac{2\bi}{\pi}(\gamma-\log 2 - \frac{3}{2}+\frac{\pi}{2}\alpha)\epsilon + \frac{2\bi}{\pi} \epsilon\log\epsilon,\\
   \label{eq:alpha}
   \alpha =& \frac{32}{\pi}<(2\rId-{\cal P}^{(e)})^{-1}\{\sqrt{m}(-1)^{m}C_0(2\pi m)\},\{\sqrt{m}(-1)^{m}C_0(2\pi m)\}>_{\ell^2},
 \end{align}
 and we recall that ${\cal P}^{(e)}$ is defined in (\ref{eq:Pe}).
 The corresponding solutions to (\ref{eq:op1}-\ref{eq:op3}) are
 \begin{align}
   b_0=&1,\\
   \{a_{2n-1}\}=& \{0\},\\
   \{a_{2n}\}=&(e^{\bi k l}-1)(2{\rm Id}-{\cal A}^{e})^{-1}\{c_{0,2n}\},
 \end{align}
 \begin{proof}
   By Lemma~\ref{lem:cmn} and by
\begin{align}
  ||(2{\rm Id} - {\cal A}^{(e)})^{-1} - (2{\rm Id}- {\cal P}^{(e)})^{-1}|| ={\cal O}(\epsilon^2\log\epsilon),
\end{align}
equation (\ref{eq:gov:k0:e}) becomes:
 \begin{align*}
  &2(e^{\bi k l}+1) = (e^{\bi k l} - 1)\Delta(\epsilon)+ {\cal O}(\epsilon^3\log\epsilon),
 \end{align*}
 which is equivalent to 
 \[
   e^{\bi k l} + 1 = -\frac{\Delta(\epsilon)}{1-\Delta(\epsilon)/2} +
   {\cal O}(\epsilon^3\log\epsilon).
 \]

 As the right-hand side approaches $0$ as $\epsilon\to 0$, we see that the
 resonance frequencies must satisfy: for some $m=1,\cdots$,
 \[
   \delta_{m,e}:= k l - k_{m,e} = o(1).
 \]
 Thus,
 \begin{align*}
   \epsilon - \epsilon_{m,e} =& h\delta_{m,e},\\ 
   \epsilon\log\epsilon - \epsilon_{m,e}\log\epsilon_{m,e}  =& h\delta_{m,e} \log\epsilon + h^2\delta_{m,e} + {\cal O}(h^3\delta_{m,e}^2k_{m,e}^{-1}),
 \end{align*}
 as $h \to 0^+$. Therefore, we have
 \[
   e^{\bi \delta_{m,e}} - 1 = \frac{\Delta(\epsilon)}{1-\Delta(\epsilon)/2} + {\cal O}(\epsilon^3\log\epsilon),
 \]
 so that by Taylor's expansion of $\log(1+x)$ and $1/(1-x)$ at $x=0$,
 \begin{align*}
   \delta_{m,e} =& -\bi\log\left[  1 +\frac{\Delta(\epsilon)}{1-\Delta(\epsilon)/2} + {\cal O}(\epsilon^3\log\epsilon)\right]\\
   =&-\bi\left[ \frac{\Delta(\epsilon)}{1-\Delta(\epsilon)/2} -\frac{\Delta^2(\epsilon)}{2(1-\Delta(\epsilon)/2)^2} + \frac{\Delta^3(\epsilon)}{3(1-\Delta(\epsilon)/2)^3} \right] + {\cal O}(\epsilon^3\log\epsilon)\\
   =&-\bi\left[\Delta(\epsilon) + \frac{1}{12}\Delta^3(\epsilon) \right] + {\cal O}(\epsilon^3\log\epsilon).
 \end{align*}
 Thus, $\delta_{m,e}\eqsim \frac{2}{\pi}\epsilon_{m,e}\log\epsilon_{m,e}$ implies
 \begin{align*}
   &\epsilon\log\epsilon - \epsilon_{m,e}\log\epsilon_{m,e} \\
   =&k_{m,e}^{-1}\epsilon_{m,e}\delta_{m,e}\log(\epsilon_{m,e}) + k_{m,e}^{-1}\epsilon_{m,e}\delta_{m,e}\log(1+\delta_{m,e}/k_{m,e}) + \epsilon_{m,e}\log(1+\delta_{m,e}/k_{m,e})\\
   =&k_{m,e}^{-1}\epsilon_{m,e}\delta_{m,e}\log(\epsilon_{m,e}) + k_{m,e}^{-1}\delta_{m,e}\epsilon_{m,e} + \frac{1}{2}k_{m,e}^{-2}\delta_{m,e}^2\epsilon_{m,e} + {\cal O}(\epsilon_{m,e}^3\log\epsilon_{m,e}).
 \end{align*}
 Based on the definition of $\Delta$, we get
 \begin{align*}
   \Delta(\epsilon)-\Delta(\epsilon_{m,e}) =& k_{m,e}^{-1}\delta_{m,e}\left[  \epsilon_{m,e} + \frac{2\bi}{\pi}(\gamma-\log 2 - \frac{3}{2}+\frac{\pi}{2}\alpha)\epsilon_{m,e}\right] \\
   &+ \frac{2\bi}{\pi} \left[  k_{m,e}^{-1}\epsilon_{m,e}\delta_{m,e}\log(\epsilon_{m,e}) + k_{m,e}^{-1}\delta_{m,e}\epsilon_{m,e} + \frac{1}{2}k_{m,e}^{-2}\delta_{m,e}^2\epsilon_{m,e}\right] + {\cal O}(\epsilon_{m,e}^3\log\epsilon_{m,e})\\
   =&k_{m,e}^{-1}\delta_{m,e}\Delta(\epsilon_{m,e}) + \frac{2\bi}{\pi} \left[  k_{m,e}^{-1}\delta_{m,e}\epsilon_{m,e} + \frac{1}{2}k_{m,e}^{-2}\delta_{m,e}^2\epsilon_{m,e}\right] + {\cal O}(\epsilon_{m,e}^3\log\epsilon_{m,e})\\
   =&{\cal O}(\epsilon_{m,e}^2\log^2\epsilon_{m,e}).
 \end{align*}
 so that
 \begin{align*}
\Delta^3(\epsilon)-\Delta^3(\epsilon_{m,e}) = (\Delta(\epsilon)-\Delta(\epsilon_{m,e}))(\Delta^2(\epsilon)+\Delta(\epsilon)\Delta(\epsilon_{m,e})+\Delta^2(\epsilon_{m,e})) = {\cal O}(\epsilon_{m,e}^4\log^4\epsilon_{m,e}).
 \end{align*}
 Therefore,
 \begin{align*}
   \delta_{m,e} =&-\bi \Delta(\epsilon_{m,e}) -\bi k_{m,e}^{-1}\delta_{m,e}\Delta(\epsilon_{m,e}) + \frac{2}{\pi} \left[  k_{m,e}^{-1}\delta_{m,e}\epsilon_{m,e} + \frac{1}{2}k_{m,e}^{-2}\delta_{m,e}^2\epsilon_{m,e}\right]\\
   &-\frac{\bi}{12}\Delta^3(\epsilon_{m,e}) + {\cal O}(\epsilon^3_{m,e}\log\epsilon_{m,e}),
 \end{align*}
 which is equivalent to 
 \[
   A \delta_{m,e}^2 + B \delta_{m,e} + C = 0,
 \]
where 
\begin{align*}
  A =& \pi^{-1}k_{m,e}^{-2}\epsilon_{m,e}={\cal O}(\epsilon_{m,e}),\\
  B =& \frac{2}{\pi} k_{m,e}^{-1}\epsilon_{m,e} - \bi k_{m,e}^{-1}\Delta(\epsilon_{m,e}) - 1 \eqsim -1,\\
  C = &-\bi \Delta(\epsilon_{m,e})-\frac{\bi}{12}\Delta^3(\epsilon_{m,e}) + {\cal O}(\epsilon^3_{m,e}\log\epsilon_{m,e})={\cal O}(\epsilon_{m,e}\log\epsilon_{m,e}).
\end{align*}
Solving this quadratic equation, 
\begin{align*}
  \delta_{m,e} =& -\frac{2C}{B+\sqrt{B^2-4AC}} = -\frac{2C}{B}\frac{1}{1+\sqrt{1-\frac{4AC}{B^2}}}\\
  =&-\frac{2C}{B}\left[ \frac{1}{2} + \frac{AC}{2B^2} \right] + {\cal O}(\epsilon_{m,e}^5\log^3\epsilon_{m,e})\\
  =&-\frac{C}{B} - \frac{AC^2}{B^3} + {\cal O}(\epsilon_{m,e}^5\log^3\epsilon_{m,e})\\
  =&-\frac{-\bi \Delta(\epsilon_{m,e})-\frac{\bi}{12}\Delta^3(\epsilon_{m,e})}{\frac{2}{\pi} k_{m,e}^{-1}\epsilon_{m,e} - \bi k_{m,e}^{-1}\Delta(\epsilon_{m,e})-1} -\frac{\pi^{-1}k_{m,e}^{-2}\epsilon_{m,e}(-\bi \Delta(\epsilon_{m,e}))^2}{(\frac{2}{\pi} k_{m,e}^{-1}\epsilon_{m,e} - \bi k_{m,e}^{-1}\Delta(\epsilon_{m,e})-1)^3} + {\cal O}(\epsilon_{m,e}^3\log\epsilon_{m,e})\\
  =&-\bi \Delta(\epsilon_{m,e})\left[1 + \left(  \frac{2}{\pi} k_{m,e}^{-1}\epsilon_{m,e} - \bi k_{m,e}^{-1}\Delta(\epsilon_{m,e})\right) + \left(  \frac{2}{\pi} k_{m,e}^{-1}\epsilon_{m,e} - \bi k_{m,e}^{-1}\Delta(\epsilon_{m,e})\right)^2  \right] \\
  &-\frac{\bi}{12}\Delta^3(\epsilon_{m,e})-\pi^{-1}k_{m,e}^{-2}\epsilon_{m,e}\Delta^2(\epsilon_{m,e}) + {\cal O}(\epsilon_{m,e}^3\log\epsilon_{m,e})\\
  =&-\bi \left[ 1 + \frac{2}{\pi} h  \right]\Delta(\epsilon_{m,e}) - \left[k_{m,e}^{-1}-5\pi^{-1}k_{m,e}^{-1}h\right]\Delta^2(\epsilon_{m,e}) + \bi\left[  k_{m,e}^{-2}-\frac{1}{12} \right]\Delta^3(\epsilon_{m,e}) \\
  &+ {\cal O}(\epsilon_{m,e}^3\log\epsilon_{m,e}).
\end{align*}
   Consequently, we see that the resonance frequency $kl$, if solving
   (\ref{eq:gov:k0:e}), must asymptotically behave as (\ref{eq:kl:e}) for
   $\epsilon_{m,e}\ll 1$. As for the existence of such solutions, one just notices
   that when $kl$ lies in $D_h=\{k\in\mathbb{C}:|kl-k_{m,e}|\leq
   h^{1/2}\}\subset {\cal S}$, then on the boundary of this disk
   \begin{align*}
     &\Bigg|2(e^{\bi k l} + 1) - (e^{\bi k l}-1)\left[  c_{00} + <({\rm Id}-{\cal A}^{(e)})^{-1}\{c_{0,2m}\},\{c_{2m,0}\}>_{\ell^2}\right] - 2\bi(kl-k_{m,e})\Bigg|\\
     =& {\cal O}(h)\leq 2\sqrt{h}=|2\bi(kl-k_{m,e})|.
   \end{align*}
   Rouch\'e's theorem indicates that there exists a unique solution to
   (\ref{eq:gov:k0:e}) in $D_h$.
  \end{proof}
\end{mytheorem}
\begin{myremark}
  In (\ref{eq:kl:e}), $\Delta^3(\epsilon_{m,e})$ contains terms greater than the
  error term ${\cal O}(\epsilon_{m,e}^3\log\epsilon_{m,e})$; we keep it here to
  make the expansion more compact and easier to evaluate. On the other hand, by (\ref{eq:kl:e}), 
  \[
    kl = k_{m,e} - \bi \left[ 1 + \frac{2}{\pi}h \right]\Delta(\epsilon_{m,e}) +
    {\cal O}(h^2\log^2 h),
  \]
  coincides with the result in Proposition 4.5 of \cite{linzha17}. By retaining
  leading behaviors of $c_{ij}$ here, we obtain an asymptotic formula of
  accuracy ${\cal O}(h^3\log h)$ that is much more accurate.
\end{myremark}
\begin{myremark}
  Like \cite{linzha17}, our asymptotic formula contains an undetermined constant
  $\alpha$ as well. In \cite{holsch19}, the authors use a method of matched
  asymptotic expansions to solve the single-slit scattering problem when a
  normal incident wave is specified and exactly describes the leading behavior
  of a real frequency (see Eq.~(35) therein) at which transmission efficiency
  reaches a peak. In fact, such real frequencies are exactly real parts of
  resonance frequencies. By comparing their formula and the real part of
  (\ref{eq:kl:e}), we easily conclude that
  \[
    \alpha=1/\pi-2/\pi\log(\pi/2).
  \]
\end{myremark}

\subsection{Odd mode}
Now, we consider the odd mode $u^o$. Since the theory is essentially the same as
the even case, we show briefly the results. $u^o$ solves
\begin{align}
  \label{eq:helm:o}
  \Delta u + k^2 u &= 0,\quad{\rm on}\quad \Omega_h,\\
  \label{eq:cond:o1}
  \partial_{\bm \nu }u &= 0,\quad{\rm on}\quad \Gamma_h\backslash \{(x_1,x_2):x_2=-l/2, |x_1|\leq h/2\},\\
  \label{eq:cond:o2}
  u &= 0,\quad{\rm on}\quad  \{(x_1,x_2):x_2=-l/2, |x_1|<h/2\}.
\end{align}
In $B_h$, we could represent $u$ as the following form,
\begin{equation}
  u(x) = \sum_{n=0}^{+\infty} b_{n}\phi_n(x_1)[e^{\bi s_n(x_2+l)}-e^{-\bi s_n x_2}].
\end{equation}
Thus, for $|x_1|\leq h/2$, we get
\[
  u_{x_2}(x_1,0) = \sum_{n=0}^{+\infty}\bi s_n b_n[e^{\bi s_nl}+1]\phi_n(x_1),
\]
so that
\[
  \hat{f}(\xi) =\sqrt{\frac{2}{h}}\sum_{n=0}^{+\infty}\bi s_{n}[e^{\bi s_{n} l} + 1]\frac{2\xi \sin(\xi h/2 + n\pi/2)}{\xi^2-\frac{\pi^2 n^2}{h^2}}b_{n}.
\]
Thus, we get
\[
  u(x_1,0) = \sum_{n=0}^{+\infty} b_n\phi_n(x_1)[e^{\bi s_n l}-1]
  =\sum_{m=0}^{+\infty} b_m[e^{\bi s_m l}+1]\psi_m(x_1).
\]
Consequently, we get the following linear system:
\begin{align}
  \label{eq:op1o}
  2b_0(e^{\bi s_0 l}-1) =& b_0(e^{\bi s_0 l} + 1)c_{00}^{(o)} + <\{a_{2m}\},\{c_{2m,0}^{(o)}\}>_{\ell^2}\\
  \label{eq:op2o}
  2\{a_{2n}\} =& b_0(e^{\bi s_0 l} + 1)\{c_{0,2n}^{(o)}\} + {\cal A}^{(e)}\{a_{2n}\},\\
  \label{eq:op3o}
  2\{a_{2n-1}\} =& {\cal A}^{(o)}\{a_{2n-1}\},
\end{align}
where we recall that $a_n=\sqrt{n}b_n$, and $c_{00}^{(o)} = \frac{4s_0h d_{00}}{\pi}$ and for $m,n\geq 1$,
\begin{align*}
c_{mn}^{(o)} =&\frac{4s_mh}{\pi}\frac{\sqrt{n}(e^{\bi s_m l}+1)}{\sqrt{m}(e^{\bi s_n l}-1)}d_{mn},\\
  c_{m0}^{(o)} =& \frac{4s_mh}{\pi}\frac{(e^{\bi s_m l}+1)}{\sqrt{m}}d_{m0},\\
  c_{0n}^{(o)} =& \frac{4s_0h\sqrt{n}}{\pi}\frac{d_{0n}}{e^{\bi s_n l}-1},
\end{align*}
and we have used $c_{mn}\equiv 0$ when $m+n\nmid 2$. Based on Lemma~\ref{lem:dmn}, we have the following lemma.
\begin{mylemma}
  \label{lem:cmno}
  For $\epsilon\ll 1$, the sequence $\{c_{mn}^{(o)}\}_{m,n=0}^{\infty}$ asymptotically behaves as:
  for $k\in\mathbb{Z}$,
  \begin{align}
    \label{eq:cmno}
    c_{mn}^{(o)} = \left\{
    \begin{array}{ll}
      0, & m+n\nmid 2;\\
      \epsilon + \frac{2\bi}{\pi}(\gamma-\log 2 - \frac{3}{2})\epsilon + \frac{2\bi}{\pi}\epsilon \log\epsilon + {\cal O}(\epsilon^3\log\epsilon) & m=n=0;\\
-4\sqrt{m}(-1)^{m/2}C_0(\pi m) + m^{-3/2}{\cal O}(\epsilon^2\log\epsilon) & n=0<m\mid 2;\\
-\frac{4\sqrt{n}}{\pi}\epsilon\bi(-1)^{n/2}C_0(\pi n) + n^{-3/2}{\cal O}(\epsilon^3\log\epsilon) & m = 0<n\mid 2;\\
p_{mn}^{(e)} -\delta_{mn} + \frac{\sqrt{mn}(\log m -\log n)}{m^2-n^2}{\cal O}(\epsilon^2) & 0<m,n\mid 2;\\
p_{mn}^{(o)} -\delta_{mn} + \frac{\sqrt{mn}(\log m -\log n)}{m^2-n^2}{\cal O}(\epsilon^2\log\epsilon) & 0<m,n\nmid 2.
    \end{array}
    \right.
  \end{align}
\end{mylemma}
Consequently, we obtain the second theorem.
\begin{mytheorem}
  \label{thm:oddres}
 For $h \ll 1$, the system (\ref{eq:op1o}-\ref{eq:op3o}) has a nonzero solution if and only if
 $k$ solves
 \begin{align}
   \label{eq:gov:k0:o}
   2(e^{\bi k l} - 1) = (e^{\bi k l}+1)\left[  c_{00}^{(o)} + <(2{\rm Id}-{\cal A}^{(e)})^{-1}\{c_{0,2m}^{(o)}\},\{c_{2m,0}^{(o)}\}>_{\ell^2}\right].
 \end{align}
 In fact, the nonzero solutions (resonance frequencies) to (\ref{eq:gov:k0:o})
 are
 \begin{align}
   \label{eq:kl:o}
   kl =& k_{m,o}-\bi \left[ 1 + \frac{2}{\pi} h  \right]\Delta(\epsilon_{m,o}) - \left[k_{m,o}^{-1}-5\pi^{-1}k_{m,o}^{-1}h\right]\Delta^2(\epsilon_{m,o}) + \bi\left[  k_{m,o}^{-2}-\frac{1}{12}\right]\Delta^3(\epsilon_{m,o}) \nonumber\\
   &+ {\cal O}(\epsilon_{m,o}^3\log\epsilon_{m,o}),\quad m=1,2,\cdots,
 \end{align}
 where $k_{m,o}=2m\pi$ are the Fabry-P\'erot frequencies and $\epsilon_{m,o}
 =k_{m,o} h\ll 1$. The corresponding solutions to (\ref{eq:op1o}-\ref{eq:op3o})
 are
 \begin{align}
   b_0=&1,\\
   \{a_{2n-1}\}=& \{0\},\\
   \{a_{2n}\}=&(e^{\bi k l}+1)(2{\rm Id}-{\cal A}^{e})^{-1}\{c_{0,2n}^{(o)}\},
 \end{align}
 \begin{proof}
   By Lemma~\ref{lem:cmn} and by
\begin{align}
  ||(2{\rm Id} - {\cal A}^{(e)})^{-1} - (2{\rm Id}- {\cal P}^{(e)})^{-1}|| ={\cal O}(\epsilon^2\log\epsilon),
\end{align}
equation (\ref{eq:gov:k0:o}) becomes:
 \begin{align*}
  &2(e^{\bi k l}-1) = (e^{\bi k l} + 1)\Delta(\epsilon)+ {\cal O}(\epsilon^3\log\epsilon),
 \end{align*}
 which is equivalent to 
 \[
   e^{\bi k l} - 1 = \frac{\Delta(\epsilon)}{1-\Delta(\epsilon)/2} +
   {\cal O}(\epsilon^3\log\epsilon).
 \]

 As the right-hand side approaches $0$ as $\epsilon\to 0$, we see that the
 resonance frequencies must satisfy: for some $m=1,\cdots$,
 \[
   \delta_{m,o}:= k l - k_{m,o} = o(1).
 \]
The proof follows from similar arguments of Theorem~\ref{thm:evenres}. We omit the details.
  \end{proof}
\end{mytheorem}

\section{Multiple slits}
\begin{figure}[!ht]
  \centering
  \includegraphics[width=0.5\textwidth]{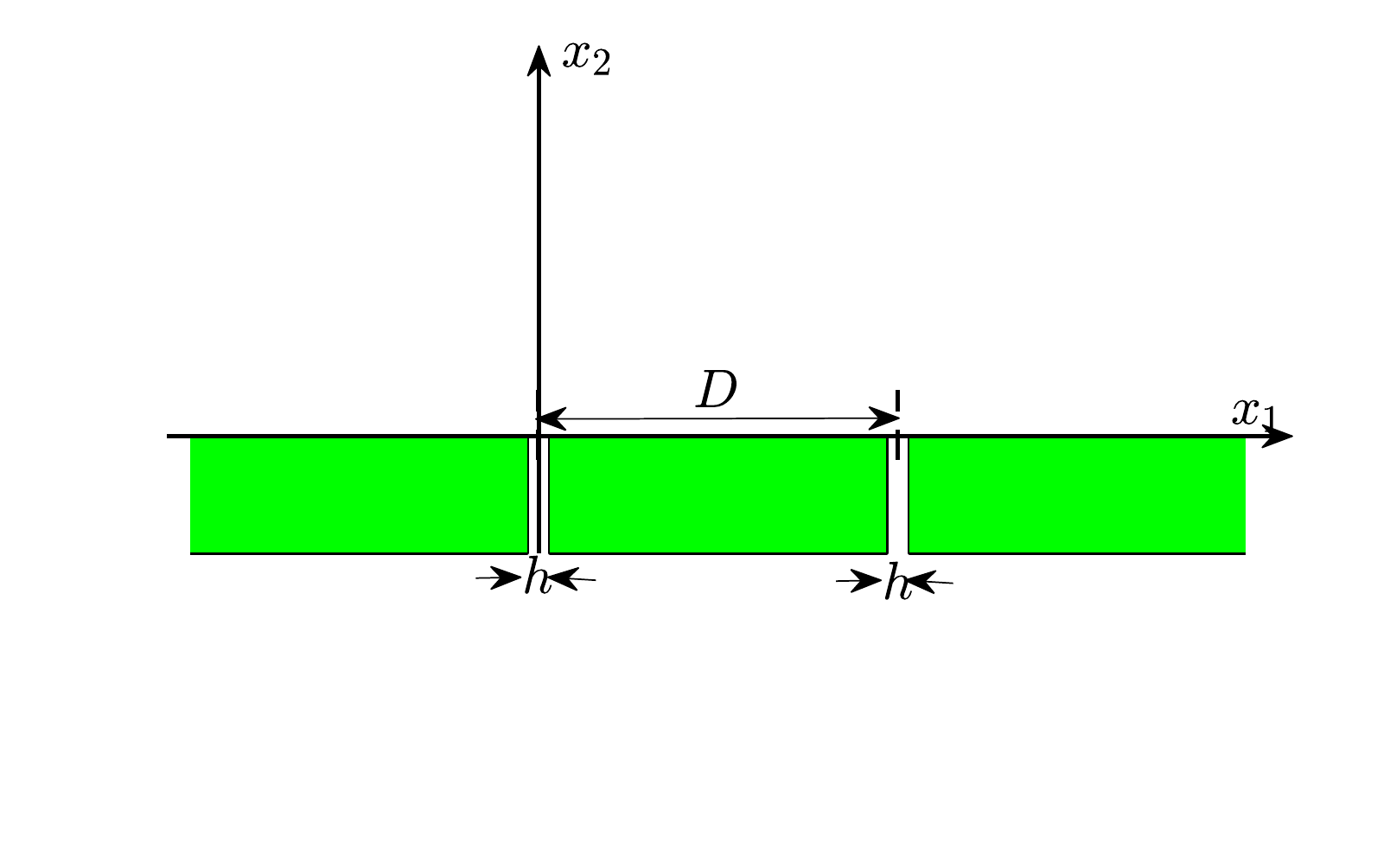}
  \caption{A perfectly conducting slab with two slits of width $h$ spaced by $D$.}
  \label{fig:model2}
\end{figure}

We study the resonance frequencies for a slab with two slits first. As we shall
see, the formula of resonance frequencies of a two-slit slab can be easily
extended to a slab with any finite number of slits. Suppose the slab has two
slits of the same width $h$ spaced by $D$ independent of $h$, as illustrated in
Figure~\ref{fig:model2}.

We consider even modes first, i.e. $u$ satisfies
$u(x_1,x_2+l/2)=u(x_1,-x_2+l/2)$. Thus, we have inside the two slits,
\begin{align*}
  u(x) =& \sum_{n=0}^{\infty}b_n\phi_n(x_1)[e^{\bi s_n(x_2+l)} + e^{-\bi s_n x_2}], |x_1|<h/2;\\
  u(x) =& \sum_{n=0}^{\infty}b_n'\phi_n(x_1-D)[e^{\bi s_n(x_2+l)} + e^{-\bi s_n x_2}], |x_1-D|<h/2.
\end{align*}
Then, we get
\begin{align*}
  \hat{f}(\xi) =& \int_{-\infty}^{+\infty}u_{x_2}(x_1,0)e^{\bi \xi x}d\xi\\
  =&\sqrt{\frac{2}{h}}\sum_{n=0}^{+\infty}\bi s_{n}[e^{\bi s_{n} l} - 1]\frac{2\xi \sin(\xi h/2 + n\pi/2)}{\xi^2-\frac{\pi^2 n^2}{h^2}}b_{n}\\
  &+\sqrt{\frac{2}{h}}\sum_{n=0}^{+\infty}\bi s_{n}[e^{\bi s_{n} l} - 1]\frac{2\xi \sin(\xi h/2 + n\pi/2)}{\xi^2-\frac{\pi^2 n^2}{h^2}}b_{n}'e^{\bi \xi D}.
\end{align*}
Thus, we obtain 
\begin{align*}
  \sum_{n=0}^{+\infty}\phi_n(x_1)b_n[e^{\bi s_n l}+1] =& \sum_{m=0}^{+\infty}(e^{\bi s_ml}-1)\left[  \psi_m(x_1)b_m + \psi_m(x_1-D)b_m'\right], |x_1|<h/2,\\
  \sum_{n=0}^{+\infty}\phi_n(x_1-D)b_n'[e^{\bi s_n l}+1] =& \sum_{m=0}^{+\infty}(e^{\bi s_ml}-1)\left[  \psi_m(x_1)b_m + \psi_m(x_1-D)b_m'\right],|x_1-D|<h/2.
\end{align*}
Taking inner product with $\phi_n$, we obtain the following linear system of
infinite dimensions,
\begin{align}
  \label{eq:op1:2holes}
  2b_0(e^{\bi s_0 l}+1) =& b_0(e^{\bi s_0 l} - 1)c_{00} + b_0'(e^{\bi s_0 l} - 1)c_{00}(D)\nonumber\\
  &+ <a_{2m},c_{2m,0}>_{\ell^2} + <a_{m}',c_{m,0}(D)>_{\ell^2}\\
  \label{eq:op2:2holes}
  2b_0'(e^{\bi s_0 l}+1) =& b_0'(e^{\bi s_0 l} - 1)c_{00} + b_0(e^{\bi s_0 l} - 1)c_{00}(-D)\nonumber\\
  &+ <a_{2m}',c_{2m,0}>_{\ell^2} + <a_{m},c_{m,0}(-D)>_{\ell^2}\\
  \label{eq:op1D}
  2\left[
  \begin{array}{c}
    \{a_{2n}\}\\
    \{a_{2n-1}\}
  \end{array}
    \right] =&b_0(e^{\bi s_0 l} - 1)\left[
  \begin{array}{c}
    \{c_{0,2n}\}\\
    \{0\}
  \end{array}
    \right]+ b_0'(e^{\bi s_0 l} - 1)\left[
  \begin{array}{c}
    \{c_{0,2n}(D)\}\\
    \{c_{0,2n-1}(D)\}
  \end{array}
    \right] \nonumber\\
  &+ \left[
  \begin{array}{cc}
    {\cal A}^{(e)} & 0\\
    0 & {\cal A}^{(o)}
  \end{array}
    \right]\left[
  \begin{array}{c}
    \{a_{2n}\}\\
    \{a_{2n-1}\}
  \end{array}
    \right] + {\cal A}(D)\{a_{n}'\}\\
  \label{eq:op2-D}
2\left[
  \begin{array}{c}
    \{a_{2n}'\}\\
    \{a_{2n-1}'\}
  \end{array}
    \right] =&b_0'(e^{\bi s_0 l} - 1)\left[
  \begin{array}{c}
    \{c_{0,2n}\}\\
    \{0\}
  \end{array}
    \right]+ b_0(e^{\bi s_0 l} - 1)\left[
  \begin{array}{c}
    \{c_{0,2n}(-D)\}\\
    \{c_{0,2n-1}(-D)\}
  \end{array}
    \right] \nonumber\\
  &+ \left[
  \begin{array}{cc}
    {\cal A}^{(e)} & 0\\
    0 & {\cal A}^{(o)}
  \end{array}
    \right]\left[
  \begin{array}{c}
    \{a_{2n}'\}\\
    \{a_{2n-1}'\}
  \end{array}
    \right] + {\cal A}(-D)\{a_{n}\}.
  \end{align}
In the above, we have defined two operators ${\cal A}(\pm D):\ell^2\to\ell^2$ such that for any $\{f_j\}_{j=1}^{\infty}\in\ell^2$,
\begin{align*}
  {\cal A}(\pm D)\{f_j\} =& \left[
                            \begin{array}{c}
                              \{\sum_{j=1}^{\infty}c_{2i,j}(\pm D)f_{j}\}_{i=1}^{\infty}\\
                              \{\sum_{j=1}^{\infty}c_{2i-1,j}(\pm D)f_{j}\}_{i=1}^{\infty}
                              \end{array}
                            \right]\in\ell^2,
\end{align*}
and $c_{00}(\pm D) = \frac{4s_0h d_{00}(\pm D)}{\pi}$ and for $m,n\geq 1$,
\begin{align*}
c_{mn}(\pm D) =&\frac{4s_mh}{\pi}\frac{\sqrt{n}(e^{\bi s_m l}-1)}{\sqrt{m}(e^{\bi s_n l}+1)}d_{mn}(\pm D),\\
  c_{m0}(\pm D) =& \frac{4s_mh}{\pi}\frac{(e^{\bi s_m l}-1)}{\sqrt{m}}d_{m0}(\pm D),\\
  c_{0n}(\pm D) =& \frac{4s_0h\sqrt{n}}{\pi}\frac{d_{0n}(\pm D)}{e^{\bi s_n l}+1},
\end{align*}
where
\begin{equation}
  d_{mn}(\pm D) =\int_{-\infty}^{+\infty}\frac{1}{\mu h^2}\frac{\xi \sin(\xi h/2 + m\pi/2)}{\xi^2-\frac{\pi^2 m^2}{h^2}}\frac{\xi \sin(\xi h/2 + n\pi/2)}{\xi^2-\frac{\pi^2 n^2}{h^2}}e^{\pm \bi \xi D} d\xi.
\end{equation}
The boundedness of ${\cal A}^{(l)}(\pm D):\ell^2\to \ell^2$ can be seen from the following
lemma.
\begin{mylemma}
  \label{lem:dmn:D}
  As $h\to 0^+$, the sequence $\{d_{mn}(\pm D)\}_{m,n=0}^{\infty}$
  asymptotically behaves as following:
  \begin{align}
    \label{eq:dmn:2holes}
    d_{mn}(\pm D) = \left\{
    \begin{array}{ll}
      \frac{\pi}{4}H_0^{(1)}(\pm kD) + {\cal O}(\epsilon^2) & m=n=0;\\
      m^{-2}{\cal O}(\epsilon^2) & n=0, m>0;\\
      n^{-2}{\cal O}(\epsilon^2) &  m=0, n>0;\\
      m^{-2}n^{-2}{\cal O}(\epsilon^2) & m, n>0,
    \end{array}
    \right.
  \end{align}
  where $H_0^{(1)}$ is the first kind Hankel function of order $0$, and the
  invisible constants in the ${\cal O}$-notations are independent of $m$ and
  $n$. Thus, $\{c_{mn}(\pm D)\}_{m,n=0}^{\infty}$ asymptotically behaves as
  following:
  \begin{align}
    \label{eq:cmn:2holes}
    c_{mn}(\pm D) = \left\{
    \begin{array}{ll}
      H_0^{(1)}(\pm kD)\epsilon  + {\cal O}(\epsilon^3) & m=n=0;\\
      m^{-3/2}{\cal O}(\epsilon^2) & n=0, 0<m;\\
      n^{-3/2}{\cal O}(\epsilon^3) &  m=0, 0<n;\\
      m^{-3/2}n^{-3/2}{\cal O}(\epsilon^2) & m, n>0.
    \end{array}
    \right.
  \end{align}
    \begin{proof}
    Here, we prove only the case $m=n=0$ as the other cases are much easier to
    justify.  We have
    \begin{align*}
      d_{00}(D) =&\int_{-\infty}^{+\infty}\frac{1}{\sqrt{k^2-\xi^2}}\frac{\sin^2(\xi h/2)}{\xi^2h^2}e^{\bi \xi D} d\xi\\ 
      =&\epsilon^{-2}\int_{+\infty\bi\to 0\to 1}\frac{1}{\sqrt{1-\xi^2}}\frac{1-\cos(\xi\epsilon)}{\xi^2}e^{\bi \xi k D}d\xi=:\epsilon^{-2} I(\epsilon;D).
    \end{align*}
    Clearly,
    \begin{align*}
      I\dd(\epsilon;D) =&\int_{+\infty\bi\to 0\to 1}\frac{1}{\sqrt{1-\xi^2}}\cos(\xi\epsilon)e^{\bi \xi k D}d\xi\\
      =&\frac{\pi}{4}\left[  H_0^{(1)}(kD+\epsilon) + H_0^{(1)}(kD-\epsilon)\right]\\
      =&\frac{\pi}{2}H_0^{(1)}(kD) + {\cal O}(\epsilon^2).
    \end{align*}
    Consequently, we get 
    \[
      d_{00}(D) = \frac{\pi}{4}H_0^{(1)}(kD) + {\cal O}(\epsilon^2).
    \]
  \end{proof}
  
\end{mylemma}
As a corollary, we obtain the following property of ${\cal A}(\pm D)$.
\begin{mycorollary}
 For $\epsilon\ll 1$, ${\cal A}(\pm D)$ is bounded from $\ell^2$ to $\ell^2$
 with   
 \[
   ||{\cal A}(\pm D)|| = {\cal O}(\epsilon^2).
 \]
\end{mycorollary}
Thus, Eqs.~(\ref{eq:op1D}-\ref{eq:op2-D}) indicate that
\begin{align}
  \left[
  \begin{array}{c}
    \{a_{2n}\}\\
    \{a_{2n-1}\}\\
    \{a_{2n}'\}\\
    \{a_{2n-1}'\}
  \end{array}
  \right] =& (e^{\bi s_0 l}-1)\epsilon\left[
  \begin{array}{c}
    b_0(2{\rm Id}-{\cal P}^{(e)})^{-1}\{-\frac{4\sqrt{2n}}{\pi}\bi(-1)^n C_0(2n\pi)\}\\
    0\\
    b_0'(2{\rm Id}-{\cal P}^{(e)})^{-1}\{-\frac{4\sqrt{2n}}{\pi}\bi(-1)^n C_0(2n\pi)\}\\
    0
  \end{array}
  \right]\nonumber\\
  &+(e^{\bi s_0 l}-1)b_0{\cal O}(\epsilon^3\log\epsilon)\{n^{-3/2}\log n \}\nonumber\\
  &+(e^{\bi s_0 l}-1)b_0'{\cal O}(\epsilon^3\log\epsilon)\{n^{-3/2}\log n \}.
\end{align}
Consequently, (\ref{eq:op1:2holes}) and (\ref{eq:op2:2holes}) are reduced to the
following linear system
\begin{align}
  \label{eq:eig:2holes}
  &\left(  \left[  2(e^{\bi k l}+1) - (e^{\bi k l}-1)\Delta(\epsilon)\right]{\rm Id}_2 - \epsilon (e^{\bi kl}-1)S_2(k,D)\right)\left[
                              \begin{array}{c}
                                b_0\\
                                b_0'
                              \end{array}
  \right] \nonumber\\
  =& (e^{\bi k l}-1)E_2(\epsilon,D)\left[
                              \begin{array}{c}
                                b_0\\
                                b_0'
                              \end{array}
                              \right],
\end{align}
where ${\rm Id}_2$ denotes the $2\times 2$ identity matrix, all four entries of the $2\times 2$ matrix $E_2(\epsilon,D)$ are ${\cal
  O}(\epsilon^3\log\epsilon)$ and
\[
S_2(k, D)=\left[
  \begin{array}{cc}
   0 &  H_0^{(1)}(kD) \\
    H_0^{(1)}(-kD) & 0
  \end{array}
\right].
\]
We have the following theorem.
\begin{mytheorem}
  \label{thm:evenres:2}
 For $h \ll 1$, the resonance frequencies of even modes of the two-slit slab are
 \begin{align}
   \label{eq:kl:e2}
   kl = &k_{m,e}-\bi \Delta_{1,j,m,e} - k_{m,e}^{-1}\Delta_{1,j,m,e}^2 + {\cal O}(\epsilon_{m,e}^2\log\epsilon_{m,e}),j=1,2; m=1,2,\cdots,
 \end{align}
 or 
 \begin{align}
   \label{eq:kl:e3}
   kl=&k_{m,e}-\bi \left[ 1 + \frac{2}{\pi} h  \right]\Delta_{2,j,m,e} - \left[k_{m,e}^{-1}-5\pi^{-1}k_{m,e}^{-1}h\right]\Delta_{2,j,m,e}^2 \nonumber\\
   &+ \bi\left[  k_{m,e}^{-2}-\frac{1}{12} \right]\Delta_{2,j,m,e}^3 + {\cal O}(\epsilon_{m,e}^3\log\epsilon_{m,e}), j=1,2; m=1,2,\cdots,
 \end{align}
 where  
 \begin{align*}
     \Delta_{1,j,m,e} =& \Delta(\epsilon_{m,e}) + \epsilon_{m,e}\lambda_j(S_2(k_{m,e},D))\\
     \delta_{1,j,m,e} =& -\bi \Delta_{1,j,m,e}(\epsilon_{m,e}) - k_{m,e}^{-1}\Delta_{1,j,m,e}^2, \\
     \Delta_{2,j,m,e} =& \Delta(\epsilon_{m,e}) + \epsilon_{m,e}\lambda_j(S_2(k_{m,e}+\delta_{1,j,m,e},D)),
 \end{align*}
 and $\lambda_j(S_2)$ indicates the $j$-th eigenvalue (in descending order of magnitude) of $S_2$ for $j=1,2$. 
 \begin{proof}
   Clearly, (\ref{eq:eig:2holes}) has a nonzero solution $[b_0,b_0']^{T}$ if and
   only if 
   \[
\left[  2(e^{\bi k l}+1) - (e^{\bi k l}-1)\Delta(\epsilon)\right]{\rm Id}_2 -
\epsilon (e^{\bi kl}-1)S_2(k,D) - (e^{\bi kl}-1) E_2(\epsilon,D),
   \]
   has a zero eigenvalue or zero determinant. Since $||\epsilon (e^{\bi kl}-1)S_2(k,D) - (e^{\bi kl}-1) E_2(\epsilon,D)||_2 = {\cal
     O}(\epsilon)$, the resonance frequency $k$ must satisfy
   \begin{align*}
     2(e^{\bi k l}+1) - (e^{\bi k l}-1)\Delta(\epsilon) = {\cal O}(\epsilon),
   \end{align*}
   so that 
   \[
     e^{\bi k l} + 1 = {\cal O}(\epsilon\log\epsilon).
   \]
   Thus, as in Theorem~\ref{thm:evenres}, 
   \[
     k = k_{m,e} + o(1),\quad h\to 0,
   \]
   for some $m=1,2,\cdots$. Obviously, $\epsilon\eqsim \epsilon_{m,e}$ and
   $\log\epsilon \eqsim \log \epsilon_{m,e}$ so that
   \[
     \delta_{m,e} = k-k_{m,e} \eqsim (-\bi) \left[  e^{\bi (k-k_{m,e})l} -
       1\right] ={\cal O}(\epsilon\log\epsilon) = {\cal O}(\epsilon_{m,e}\log\epsilon_{m,e}).
   \]
   Thus,
   \[
2(e^{\bi k l}+1) - (e^{\bi k l}-1)\Delta(\epsilon) - \epsilon (e^{\bi
  kl}-1)\lambda_j(S_2(k_{m,e},D))  = {\cal O}(\epsilon^2\log\epsilon), j=1,2,
   \]
   where $\lambda_j(S_2(k_{m,e},D))$ denotes the $j$-th eigenvalue of
   $S_2(k_{m,e},D)$ for $j=1,2$. By the same procedures in Theorem~\ref{thm:evenres}, we get
   \[
     \delta_{m,e} = \delta_{1,j,m,e} +{\cal
       O}(\epsilon_{m,e}^2\log\epsilon_{m,e}),
   \]
   where
   \begin{align*}
     \delta_{1,j,m,e} =& -\bi \Delta_{1,j,m,e} - k_{m,e}^{-1}\Delta_{1,j,m,e}^2, \\
     \Delta_{1,j,m,e} =& \Delta(\epsilon_{m,e}) + \epsilon_{m,e}\lambda_j(S_2(k_{m,e},D)).
   \end{align*}
   Now, we have
   \[
2(e^{\bi k l}+1) - (e^{\bi k l}-1)\Delta(\epsilon) - \epsilon (e^{\bi
  kl}-1)\lambda_j(S_2(k_{m,e}+\delta_{1,j,m,e},D))  = {\cal O}(\epsilon^3\log\epsilon),
   \]
   where we note that the $j$-th eigenvalue $\lambda_j$ of
   $S_2(k_{m,e}+\delta_{1,j,m,e},D)$ should be the one close to $\lambda_j$ with
   distance $O(\delta_{1,j,m,e})$. Then, following similar arguments in
   Theorem~\ref{thm:evenres} again, we get
   \begin{align*}
     \delta_{m,e} =&-\bi \left[ 1 + \frac{2}{\pi} h  \right]\Delta_{2,j,m,e} - \left[k_{m,e}^{-1}-5\pi^{-1}k_{m,e}^{-1}h\right]\Delta_{2,j,m,e}^2 + \bi\left[  k_{m,e}^{-2}-\frac{1}{12} \right]\Delta_{2,j,m,e}^3\\
  &+ {\cal O}(\epsilon_{m,e}^3\log\epsilon_{m,e}).
   \end{align*}

   We now prove the existence of those solutions. Since $S_2(k_m,D)$ is skew-Hermitian, one could find
   a unitary matrix $Q$, s.t., $Q^*S_2(k_m,D)Q={\rm
     diag}\{\lambda_1(S_2(k_m,D)),\lambda_2(S_2(k_m,D))\}$. Then,
   Eq.~(\ref{eq:eig:2holes}) becomes
   \begin{align*}
  &\left(  \left[  2(e^{\bi k l}+1) - (e^{\bi k l}-1)\Delta(\epsilon)\right]{\rm Id}_2 - \epsilon (e^{\bi kl}-1){\rm diag}\{\lambda_1(S_2(k_m,D)),\lambda_2(S_2(k_m,D))\}\right)\left[
                              \begin{array}{c}
                                \tilde{b}_0\\
                                \tilde{b}_0'
                              \end{array}
  \right] \nonumber\\
  =& (e^{\bi k l}-1)Q^*[E_2(\epsilon,D) + \epsilon(S_2(k,D)-S_2(k_m,D))]Q\left[
                               \begin{array}{c}
                                \tilde{b}_0\\
                                \tilde{b}_0'
                               \end{array}
     \right],
   \end{align*}
   where 
   \[
\left[
                               \begin{array}{c}
                                \tilde{b}_0\\
                                \tilde{b}_0'
                               \end{array}
     \right]
     =Q^*\left[
                               \begin{array}{c}
                                \tilde{b}_0\\
                                \tilde{b}_0'
                               \end{array}
     \right].
   \]
   Assume that $k$ lies in the disk $D_h=\{k\in\mathbb{C}:|kl-k_{m,e}|\leq
   h^{1/2}\}$. Then, on the boundary of $D_h$, all entries of
   \[
     Q^*[E_2(\epsilon,D) + \epsilon(S_2(k,D)-S_2(k_m,D))]Q,
   \]
   are ${\cal O}(h^{3/2})$,
   so that by the linearity of determinant,
   \[
     \left|{\rm Det}_1 - {\rm Det}_2\right|={\cal O}(h^2)\leq {\cal O}(h)=\left| {\rm Det}_2 \right|,
   \]
   where
   \begin{align*}
     {\rm Det}_1 = &\Bigg|  {\rm diag}\left\{  2(e^{\bi k l}+1) - (e^{\bi k l}-1)\Delta(\epsilon)- \epsilon (e^{\bi kl}-1)\lambda_j(S_2(k_m,D))\right\}_{j=1}^2 \\
     &- (e^{\bi k l}-1)Q^*[E_2(\epsilon,D) + \epsilon(S_2(k,D)-S_2(k_m,D))]Q\Bigg|,\\
     {\rm Det}_2 = &\Bigg|  {\rm diag}\left\{  2(e^{\bi k l}+1) - (e^{\bi k l}-1)\Delta(\epsilon)- \epsilon (e^{\bi kl}-1)\lambda_j(S_2(k_m,D))\right\}_{j=1}^2 \Bigg|.
   \end{align*}
   For either $j=1,2$, it is clear that on the boundary of $D_h$,
   \begin{align*}
     &\left|  2(e^{\bi k l}+1) - (e^{\bi k l}-1)\Delta(\epsilon)- \epsilon (e^{\bi kl}-1)\lambda_j(S_2(k_m,D))- 2\bi (kl - k_{m,e}) \right| \\
     =& {\cal O}(h)\leq |2\bi(kl-k_{m,e})|.
   \end{align*}
   The above two inequalities and Rouch\'e's theorem indicate that there 
   are exactly two solutions in $D_h$.
 \end{proof}
\end{mytheorem}
The following theorem characterizes resonance frequencies of odd modes, i.e.,
when the wave field $u$ satisfies $u(x_1,-x_2+l/2)=-u(x_1,x_2+l/2)$.
\begin{mytheorem}
  \label{thm:oddres:2}
 For $h \ll 1$, the resonance frequencies of odd modes of the two-slit slab are
 \begin{align}
   \label{eq:kl:o2}
   kl = &k_{m,o}-\bi \Delta_{1,j,m,o} - k_{m,o}^{-1}\Delta_{1,j,m,o}^2 + {\cal O}(\epsilon_{m,o}^2\log\epsilon_{m,o}),j=1,2; m=1,2,\cdots,
 \end{align}
 or 
 \begin{align}
   \label{eq:kl:o3}
   kl=&k_{m,o}-\bi \left[ 1 + \frac{2}{\pi} h  \right]\Delta_{2,j,m,o} - \left[k_{m,o}^{-1}-5\pi^{-1}k_{m,o}^{-1}h\right]\Delta_{2,j,m,o}^2 + \bi\left[  k_{m,o}^{-2}-\frac{1}{12} \right]\Delta_{2,j,m,o}^3 \nonumber\\
  &+ {\cal O}(\epsilon_{m,o}^3\log\epsilon_{m,o}), j=1,2; m=1,2,\cdots,
 \end{align}
 where  
 \begin{align*}
     \Delta_{1,j,m,o} =& \Delta(\epsilon_{m,o}) + \epsilon_{m,o}\lambda_j(S_2(k_{m,o},D))\\
     \delta_{1,j,m,o} =& -\bi \Delta_{1,j,m,o} - k_{m,o}^{-1}\Delta_{1,j,m,o}^2, \\
     \Delta_{2,j,m,o} =& \Delta(\epsilon_{m,o}) + \epsilon_{m,o}\lambda_j(S_2(k_{m,o}+\delta_{1,j,m,o}^1,D)),
 \end{align*}
 and $\lambda_j(S_2)$ indicates the $j$-th eigenvalue (in descending order of magnitude) of $S_2$ for $j=1,2$. 
 \begin{proof}
   The proof follows from similar arguments as in Theorem~\ref{thm:evenres:2}.
 \end{proof}
 \end{mytheorem}
 The above results can be readily extended to a slab with three or more
 slits. Specifically, suppose now the slab has $N$ slits of the same width $h$
 and thickness $l$, centered at
 \[
   (D_1,-l/2), (D_2,-l/2),\cdots, (D_{N},-l/2),
 \]
 respectively. We state our main result in the following. 
 \begin{mytheorem}
   \label{thm:res:N}
 For $h \ll 1$, the resonance frequencies of the $N$-slit slab are
 \begin{align}
   \label{eq:kl:N2}
   kl = &k_{m}-\bi \Delta_{1,j,m} - k_{m}^{-1}\Delta_{1,j,m}^2 + {\cal O}(\epsilon_{m}^2\log\epsilon_{m}),j=1,\cdots, N; m=1,2,\cdots,
 \end{align}
 or 
 \begin{align}
   \label{eq:kl:N3}
   kl=&k_{m}-\bi \left[ 1 + \frac{2}{\pi} h  \right]\Delta_{2,j,m} - \left[k_{m}^{-1}-5\pi^{-1}k_{m}^{-1}h\right]\Delta_{2,j,m}^2 + \bi\left[  k_{m}^{-2}-\frac{1}{12} \right]\Delta_{2,j,m}^3 \nonumber\\
  &+ {\cal O}(\epsilon_{m}^3\log\epsilon_{m}), j=1,\cdots, N; m=1,2,\cdots,
 \end{align}
 where $k_m=m\pi$ are the Fabry-P\'erot frequencies, $\epsilon_m=k_mh\ll 1$,
 \begin{align*}
     \Delta_{1,j,m} =& \Delta(\epsilon_{m}) + \epsilon_m\lambda_j(S_N(k_{m},\{D_j\}_{j=1}^{N}))\\
     \delta_{1,j,m} =& -\bi \Delta_{1,j,m}(\epsilon_{m}) - k_{m}^{-1}\Delta_{1,j,m}^2, \\
     \Delta_{2,j,m} =& \Delta(\epsilon_m) + \epsilon_m\lambda_j(S_N(k_{m}+\delta_{1,j,m},\{D_j\}_{j=1}^{N})),
 \end{align*}
and $\lambda_j(S_N(k,\{D_j\}_{j=1}^N))$ indicates the $j$-th eigenvalue (in descending order of magnitude) of
 \begin{align}
   &S_N(k,\{D_j\}_{j=1}^{N}) \nonumber\\
   =& \left[
  \begin{array}{ccccc}
    0 &  H_0^{(1)}(kD_{12}) & \cdots & H_0^{(1)}(kD_{1,N-1}) &  H_0^{(1)}(kD_{1N})\\
    H_0^{(1)}(kD_{21}) & 0 & \cdots & H_0^{(1)}(kD_{2,N-1}) & H_0^{(1)}(kD_{2N}) \\
    \vdots & \vdots & \vdots & \vdots & \vdots\\
    H_0^{(1)}(kD_{N1}) & H_0^{(1)}(kD_{N2})&\cdots & H_0^{(1)}(kD_{N,N-1}) & 0\\
  \end{array}
\right],
 \end{align}
with $D_{ij}=D_j-D_i$ for $i,j=1,\cdots,N$.
 \begin{proof}
   The proof is analogous to that of Theorems~\ref{thm:evenres:2} and \ref{thm:oddres:2}.
 \end{proof}
\end{mytheorem}
\begin{myremark}
  Since $S_N(k_{m},\{D_j\}_{j=1}^N)$ is skew-Hermitian, its eigenvalues are all
  pure imaginary, so that we see from (\ref{eq:kl:N2}) in
  Theorem~\ref{thm:res:N} that the imaginary part of any resonance frequency is
  asymptotically leaded by $-k_mh$, and can never attain ${\cal O}(h^2)$ however
  those slits are placed as long as $D_{ij}\gg h$ in matrix
  $S_N(k,\{D_j\}_{j=1}^{N})$. In other words, to make ${\rm Im}(k)\ll h$, at
  least one of $D_{ij}$ should be comparable to $h$, as was illustrated by
  Babadjian et al. \cite{babbontri10}, where the structure contains two slits
  spaced by ${\cal O}(h)$.
\end{myremark}
\section{Conclusion}
We have proposed a quite simple Fourier-transformation approach to study
resonances in a perfectly conducting slab with finite number of subwavelength
slits of width $h\ll 1$. Outside the slits, we Fourier transformed the governing
equation and expressed wave field in terms of field derivatives on the aperture.
Inside the slit, wave field was expressed as Fourier series in terms of a
countable basis functions with unknown Fourier coefficients. By matching field
on the aperture, we established a linear system of infinite number of equations
governing the countable Fourier coefficients. By asymptotic analysis of each
entry of the coefficient matrix, we have rigorously shown that, by removing only
a finite number of rows and columns, the resulting principle sub-matrix is
diagonally dominant so that the infinite dimensional linear system is reduced to
a finite dimensional linear system. This in turn provided a simple, asymptotic
formula of resonance frequencies of accuracy ${\cal O}(h^3\log h)$. This
asymptotic formula rigorously confirms a fact that the imaginary part of
resonance frequencies is always ${\cal O}(h)$ no matter how we place the slits
as long as they are spaced by distances independent of $h$.

As no subregion Green functions are required, we could see immediate advantages
in analyzing more complicated structures. For a slab with impedance boundary
condition, the background Green function involves Sommerfeld integrals which are
not easy to analyze. For a slab with periodic slits, our theory does not need to
evaluate the quasi-periodic Green function. Furthermore, for an ideal PEC slab
with a circular or rectangular hole, the Green function of the region of the
hole could be quite challenging to derive or analyze. Therefore, we expect that
our method could serve as an efficient approach in analyzing resonances in such
structures. We shall report the results in a future work.
\section*{Appendix}
   
To study the asymptotic behavior of $d_{mn}$ for $\epsilon=kh\ll 1$ where
$k\in{\cal S}$,
we need the following technical lemmas.
\begin{mylemma}
 \label{lem:asym:I} 
 Let $b>1$ and
 \begin{align}
   \label{eq:int:1}
   I(\epsilon;b) = \int_1^{+\infty}\frac{1-e^{-\epsilon t}}{\sqrt{1+t^2}(\epsilon^2t^2+b^2)}dt.
 \end{align}
 As $h\to 0^+$, 
 \begin{align*}
    I(\epsilon;b)=&C_0(b) - \sqrt{2}\epsilon b^{-2} + b^{-2}{\cal O}(\epsilon^2 \log\epsilon),
 \end{align*}
 where the invisible constants in the ${\cal O}$-notation are independent of
 $b$ and $k$, and we recall that $C_0(b)$ is defined in (\ref{eq:def:C0}). 
 \begin{proof}
   We make the rescaling $\epsilon = \frac{k}{|k|}(|k|h)$ so
   that one could assume $|k|=1$ in the following.
   First,
   \begin{align*}
     I(\epsilon;b) =& \int_{h}^{\infty}\frac{1-e^{-kt}}{\sqrt{h^2+t^2}(k^2t^2+b^2)}dt \\
     &= \int_{h}^{1}+\int_{1}^{\infty}\frac{1-e^{-kt}}{\sqrt{h^2+t^2}(k^2t^2+b^2)}dt\\
&=:I_1(\epsilon;b) + I_2(\epsilon;b).
   \end{align*}
   Here,
   \begin{align*}
    I_2(\epsilon;b) =& \int_{1}^{\infty}\frac{1-e^{-kt}}{\sqrt{h^2+t^2}(k^2t^2+b^2)}dt\\ 
     =&\int_{1}^{\infty}\frac{1-e^{-kt}}{t(k^2t^2+b^2)}dt - \int_1^\infty\frac{1-e^{-kt}}{(k^2t^2+b^2)}\left[ \frac{h^2}{\sqrt{h^2+t^2}t(\sqrt{h^2+t^2}+t)} \right]dt\\
     =&\int_{1}^{\infty}\frac{1-e^{-kt}}{t(k^2t^2+b^2)}dt + b^{-2}{\cal O}(h^2).
   \end{align*}
   On the other hand,
   \begin{align*}
     I_1(\epsilon;b) =& \int_{h}^{1}\frac{1-e^{-kt}}{\sqrt{h^2+t^2}(k^2t^2+b^2)}dt\\
     =&\int_{0}^{1}\frac{1-e^{-kt}-kt+\frac{k^2t^2}{2}}{t(k^2t^2+b^2)}dt - \int_{0}^{h}\frac{1-e^{-kt}-kt+\frac{k^2t^2}{2}}{t(k^2t^2+b^2)}dt \\
     &- \int_{h}^{1}\frac{1-e^{-kt}-kt+\frac{k^2t^2}{2}}{k^2t^2+b^2}\left[ \frac{h^2}{\sqrt{h^2+t^2}t(\sqrt{h^2+t^2}+t)} \right]dt  \\
     &+ \int_{h}^{1}\frac{kt-\frac{k^2t^2}{2}}{\sqrt{h^2+t^2}(k^2t^2+b^2)}dt\\
     =&\int_{0}^{1}\frac{1-e^{-kt}-kt+\frac{k^2t^2}{2}}{t(k^2t^2+b^2)}dt + b^{-2}{\cal O}(h^3) + b^{-2}{\cal O}(h^2)\\
     &+ b^{-2}\int_{h}^{1}\frac{kt-\frac{k^2t^2}{2}}{\sqrt{h^2+t^2}}dt -\int_{h}^{1}\frac{k^3t^3-\frac{k^4t^4}{2}}{t(k^2t^2+b^2)b^2}dt \\
     &+\int_{h}^{1}\frac{k^3t^3-\frac{k^4t^4}{2}}{(k^2t^2+b^2)b^2}\left[ \frac{h^2}{\sqrt{h^2+t^2}t(\sqrt{h^2+t^2}+t)} \right]dt \\
     =&\int_{0}^{1}\frac{1-e^{-kt}-kt+\frac{k^2t^2}{2}}{t(k^2t^2+b^2)}dt -\int_{0}^{1}\frac{k^3t^2-\frac{k^4t^3}{2}}{(k^2t^2+b^2)b^2}dt\\
     &+ b^{-2}\left[ k\sqrt{h^2+1} - k\sqrt{2}h -\frac{k^2}{2}\left(\frac{\sqrt{h^2+1}}{2}-\frac{\sqrt{2}}{2}h^2 - \frac{h^2}{2}\sinh^{-1}(h^{-1}) + \frac{h^2}{2}\sinh^{-1}(1)  \right) \right]\\
     &+ b^{-2}{\cal O}(h^2).
   \end{align*}
   Combining the above yields the asymptotic behavior of $I(\epsilon;b)$.
 \end{proof}
\end{mylemma}
\begin{mylemma}
  \label{lem:K+-}
  Let $b\geq b'>e$ and 
  \begin{align}
    K^{\pm}(\epsilon;b,b')=\int_{1}^{\infty}\frac{\epsilon^2 t^2(1\pm e^{-\epsilon t})}{\sqrt{1+t^2}(\epsilon^2t^2+b^2)(\epsilon^2 t^2+b'^2)}dt,
  \end{align}
  then
  \begin{align*}
    K^{-}(\epsilon;b,b')=&C_0^{-}(b,b') + \frac{\log b -\log b'}{b^2-b'^2}{\cal O}(\epsilon^2),\\
    K^{+}(\epsilon;b,b')=&C_0^{+}(b,b') + \frac{\log b -\log b'}{b^2-b'^2}{\cal O}(\epsilon^2\log\epsilon),
  \end{align*}
  where 
  the invisible constants in the big-${\cal O}$ notation are independent of $b$
and $b'$, and $C_0^{\pm}(b,b')$ are defined in (\ref{eq:def:C0-}) and
(\ref{eq:def:C0+}), respectively.
  \begin{proof}
   As in the previous lemma, we could assume $|k|=1$ and $k=k_1-\bi k_2$ with
   $k_1\geq k_2>0$. We have
   \begin{align*}
     K^-(\epsilon;b,b') =& \int_{h}^{\infty}\frac{k^2 t^2(1- e^{-k t})}{\sqrt{h^2+t^2}(k^2t^2+b^2)(k^2 t^2+b'^2)}dt\\
     =&\int_{h}^{1}+\int_{1}^{\infty}\frac{k^2 t^2(1- e^{-k t})}{\sqrt{h^2+t^2}(k^2t^2+b^2)(k^2 t^2+b'^2)}dt\\
     =&K_1^-(\epsilon;b,b') + K_2^-(\epsilon;b,b').
   \end{align*}
   Thus,
   \begin{align*}
     &\left|  K_2^-(\epsilon;b,b') - \int_{1}^{\infty}\frac{k^2 t^2(1- e^{-k t})}{t(k^2t^2+b^2)(k^2 t^2+b'^2)}dt\right| \\
     =&\left|  \int_{1}^{\infty}\frac{k^2 t^2(1- e^{-k t})}{(k^2t^2+b^2)(k^2 t^2+b'^2)}\left[ \frac{h^2}{t\sqrt{h^2+t^2}(t+\sqrt{h^2+t^2})} \right]dt\right|\\
     \leq &2 h^2\int_{1}^{\infty}\frac{1}{[(k_1^2-k_2^2)t^2+b^2][(k_1^2-k_2^2)t^2+b'^2]t}dt \\
     =&h^2\frac{(b^2\log(b'^2 + (k_1^2-k_2^2)) - b'^2\log(b^2 + (k_1^2-k_2^2)))}{b^2b'^2(b^2 - b'^2)} + h^2\frac{\log (k_1^2-k_2^2)}{b^2b'^2}\\
     =& \frac{\log b - \log b'}{b^2-b'^2}{\cal O} (h^2).
   \end{align*}
   Moreover,
   \begin{align*}
    &\left|  K_1^-(\epsilon;b,b') - 
      \int_{0}^{1}\frac{k^2 t^2(1- e^{-k t})}{t(k^2t^2+b^2)(k^2 t^2+b'^2)}dt\right| \\
    =&\left|\int_{0}^{h}\frac{k^2 t^2(1- e^{-k t})}{t(k^2t^2+b^2)(k^2 t^2+b'^2)}dt-\int_{h}^{1}\frac{k^2 t^2(1- e^{-k t})}{(k^2t^2+b^2)(k^2 t^2+b'^2)}\left[  \frac{h^2}{t\sqrt{h^2+t^2}(t+\sqrt{h^2+t^2})}\right]dt\right|\\
     =&b^{-2}b'^{-2}{\cal O}(h^2),
   \end{align*}
   which yields the desired results for $K^-$. Similarly, one obtains 
   \begin{align*}
     K^+(\epsilon;b,b') =& K^-(\epsilon;b,b') + \int_{h}^{\infty}\frac{2k^2 t^2e^{-kt}}{\sqrt{h^2+t^2}(k^2t^2+b^2)(k^2 t^2+b'^2)}dt\\
     =&K^-(\epsilon;b,b') + \int_{0}^{1}\frac{2k^2 t^2e^{-kt}}{t(k^2t^2+b^2)(k^2 t^2+b'^2)}dt + \int_0^h\frac{2k^2 t^2e^{-kt}}{t(k^2t^2+b^2)(k^2 t^2+b'^2)}dt\\
     &+ \int_{h}^{1}\frac{2k^2 t^2e^{-kt}}{(k^2t^2+b^2)(k^2 t^2+b'^2)}\frac{h^2}{t\sqrt{h^2+t^2}(t+\sqrt{h^2+t^2})}dt \\
     &+\int_{1}^{\infty}\frac{2k^2 te^{-kt}}{(k^2t^2+b^2)(k^2 t^2+b'^2)}dt+ \frac{\log b - \log b'}{b^2-b'^2}{\cal O} (h^2)\\
     =&C_0^{-}(b,b') +  \int_{0}^{\infty}\frac{2k^2 t^2e^{-kt}}{t(k^2t^2+b^2)(k^2 t^2+b'^2)}dt \\
     &+ b^{-2}b'^{-2}{\cal O}(h^2\log h)+ \frac{\log b - \log b'}{b^2-b'^2}{\cal O} (h^2)\\
     =&\int_0^{\infty}\frac{t(1+e^{-t})}{(t^2+b^2)(t^2+b'^2)}dt + \frac{\log b - \log b'}{b^2-b'^2}{\cal O} (h^2\log h),
   \end{align*}
   which concludes the proof.
  \end{proof}

\end{mylemma}
By the above two lemmas, we are ready to analyze the asymptotic behavior of
$d_{mn}$ for $\epsilon\ll 1$.

\noindent\textbf{Proof of Lemma~\ref{lem:dmn}.}  It is clear that $d_{mn}=0$ when $m+n\nmid 2$. For the case $m+n\mid 2$, we
have for $k\in\mathbb{R}^+$ that
\begin{align*}
  d_{mn}=&\int_{-\infty}^{+\infty}\frac{1}{\mu h^2}\frac{\xi \sin(\xi h/2 + m\pi/2)}{\xi^2-\frac{\pi^2 m^2}{h^2}}\frac{\xi \sin(\xi h/2 + n\pi/2)}{\xi^2-\frac{\pi^2 n^2}{h^2}} d\xi\\
  =&\frac{(-1)^{(m-n)/2}}{2}\int_{-\infty}^{\infty}\frac{\epsilon^2\xi^2(1-(-1)^m\cos(\epsilon\xi))}{\sqrt{1-\xi^2}(\epsilon^2\xi^2-\pi^2m^2)(\epsilon^2\xi^2-\pi^2n^2)}d\xi.
\end{align*}
We consider case $m=n=0$ first. Then, we have by Cauchy's theorem that
\begin{align}
  \label{eq:newdef:d00}
  d_{00}=&\frac{1}{2}\int_{-\infty}^{+\infty}\frac{1-\cos(\epsilon\xi)}{\sqrt{1-\xi^2}\epsilon^2\xi^2}d\xi =\frac{1}{2}\int_{-\infty}^{+\infty}\frac{1-e^{\bi\epsilon\xi}+\bi \epsilon\xi}{\sqrt{1-\xi^2}\epsilon^2\xi^2}d\xi =\int_{+\infty\bi\to 0\to 1}\frac{1-e^{\bi\epsilon\xi}+\bi \epsilon\xi}{\sqrt{1-\xi^2}\epsilon^2\xi^2}d\xi\nonumber\\
  =&\int_{0}^1\frac{1+\bi\epsilon\xi-e^{\bi\epsilon\xi}}{\sqrt{1-\xi^2}\epsilon^2\xi^2}d\xi + \epsilon^{-2}\bi \int_{0}^{+\infty}\frac{1-e^{-\epsilon t}-\epsilon t}{\sqrt{1+t^2}{t^2}}dt\\
  =:&I_1(\epsilon) + \epsilon^{-2}I_2(\epsilon)\nonumber.
\end{align}
Here by uniqueness principle, (\ref{eq:newdef:d00}) holomorphically extends the
definition of $d_{00}(k)$ from $\mathbb{R}^+$ to ${\cal S}$.
Now, by Taylor series of $e^{\bi \epsilon \xi}$, we get the asymptotic
expansion for $I_1$ as follows
\begin{align*}
  I_1(\epsilon)=&\frac{\pi}{4} + \frac{\bi \epsilon}{6} + {\cal O}(\epsilon^2).
\end{align*}
On the other hand,
\begin{align*}
  I_2''(\epsilon)=&-\bi\int_{0}^{+\infty}\frac{e^{- \epsilon t}}{\sqrt{1+t^2}}dt\\
  =& \bi \frac{\pi}{2} Y_0(\epsilon)-\bi \int_0^1(1-t^2)^{-1/2}\sin(\epsilon t)dt\\
  =& \bi(\log(\epsilon/2)+\gamma) - \bi\epsilon + {\cal O}(\epsilon^2\log(\epsilon)),
\end{align*}
where $\gamma$ is the Euler's constant. Thus,
\[
  I_2(\epsilon) = \bi\frac{\gamma-\log 2}{2}\epsilon^2 -\bi \frac{\epsilon^3}{6}
  + \bi\left[\frac{\epsilon^2\log\epsilon}{2}-\frac{3\epsilon^2}{4}\right] + {\cal O}(\epsilon^4\log\epsilon).
\]
Consequently, we get
\[
  d_{00}= \frac{\pi}{4} + \frac{\bi}{2}(\gamma-\log 2 - \frac{3}{2}) + \frac{\bi}{2} \log\epsilon + {\cal O}(\epsilon^2\log\epsilon).
\]
Now, consider the case when $m=0$ and $0\neq n\mid 2$. We have
\begin{align*}
  d_{n0}=d_{0n} = &(-1)^{n/2}\left[ \int_0^1 \frac{1-e^{\bi\epsilon\xi}}{\sqrt{1-\xi^2}(\epsilon^2\xi^2-\pi^2n^2)}dt +\bi\int_0^{\infty}\frac{1-e^{-\epsilon t}}{\sqrt{1+t^2}(\epsilon^2t^2+\pi^2n^2)}dt \right]\\
  =:&(-1)^{n/2}\left[I_3(\epsilon) + \bi I_4(\epsilon)\right].
\end{align*}
Clearly,
\begin{align*}
  I_3(\epsilon) = \frac{\bi \epsilon}{\pi^2n^2} + n^{-2}{\cal O}(\epsilon^2). 
\end{align*}
On the other hand, according to Lemma~\ref{lem:asym:I}, 
\begin{align*}
  I_4(\epsilon) =& \int_{0}^{1}\frac{1-e^{-\epsilon t}}{\sqrt{1+t^2}(\epsilon^2t^2+\pi^2n^2)}dt + I(\epsilon;\pi n)\\ 
  =&\frac{\epsilon}{\pi^2n^2}(\sqrt{2}-1) + n^{-2}{\cal O}(\epsilon^3) + C_0(\pi n) - \frac{\sqrt{2}\epsilon}{\pi^{2}n^{2}} + n^{-2}{\cal O}(\epsilon^2\log\epsilon)\\
  =&C_0(\pi n) -\frac{\epsilon}{\pi^2n^2}+n^{-2}{\cal O}(\epsilon^2\log\epsilon).
\end{align*}
Consequently, we get
\[
  d_{n0}=d_{0n} = \bi(-1)^{n/2}C_0(\pi n) + n^{-2}{\cal O}(\epsilon^2\log\epsilon).
\]
When $mn\neq 0$ and $m\neq n$, we have by Cauchy's theorem that
\begin{align*}
  d_{mn}=&\frac{1}{2}(-1)^{(m-n)/2}\int_{-\infty}^{\infty}\frac{\epsilon^2\xi^2(1-(-1)^me^{\bi\epsilon\xi})}{\sqrt{1-\xi^2}(\epsilon^2\xi^2-\pi^2m^2)(\epsilon^2\xi^2-\pi^2n^2)}d\xi\\
  =&(-1)^{(m-n)/2}\Big[  \int_{0}^{1}\frac{\epsilon^2\xi^2(1-(-1)^me^{\bi\epsilon\xi})}{\sqrt{1-\xi^2}(\epsilon^2\xi^2-\pi^2m^2)(\epsilon^2\xi^2-\pi^2n^2)}d\xi\\
         &+ \bi\int_{0}^{\infty}\frac{\epsilon^2t^2(1-(-1)^me^{-\epsilon t})}{\sqrt{1+t^2}(\epsilon^2t^2+\pi^2m^2)(\epsilon^2t^2+\pi^2n^2)}dt\Big]\\
  =:&(-1)^{(m-n)/2}\left[ I_5(\epsilon) + \bi I_6(\epsilon)\right].
\end{align*}
Clearly, we have
\[
  I_5(\epsilon) = \frac{\epsilon^2}{\pi^4m^2n^2}\left[  (1-(-1)^m)\frac{\pi}{4}
    -(-1)^m\frac{2}{3}\bi \epsilon\right] + m^{-2}n^{-2}{\cal O}(\epsilon^4).
\]
If $m\mid 2$, we have
\begin{align*}
  I_6(\epsilon) =& \int_{0}^{1}\frac{\epsilon^2t^2(1-e^{-\epsilon t})}{\sqrt{1+t^2}(\epsilon^2t^2+\pi^2m^2)(\epsilon^2t^2+\pi^2n^2)}dt + K^{-}(\epsilon;\pi m,\pi n)\\
  =&m^{-2}n^{-2}{\cal O}(\epsilon^{3})+C_0^{-}(\pi m,\pi n) + \frac{\log m -\log n}{m^2-n^2}{\cal O}(\epsilon^2).
\end{align*}
Consequently,
\[
  d_{mn} = \bi(-1)^{(m-n)/2}C_0^-(\pi m,\pi n) + \frac{\log m -\log n}{m^2-n^2}{\cal O}(\epsilon^2).
\]
If $m\nmid 2$, we similarly have
\[
  d_{mn} = \bi(-1)^{(m-n)/2}C_0^+(\pi m,\pi n) + \frac{\log m -\log n}{m^2-n^2}{\cal O}(\epsilon^2\log\epsilon).
\]
If $0\neq m=n \mid 2$, we have by Cauchy's theorem that 
\begin{align*}
  d_{mm}&=\frac{1}{2}\int_{-\infty}^{\infty}\frac{\epsilon^2\xi^2(1-\cos(\epsilon\xi))}{\sqrt{1-\xi^2}(\epsilon^2\xi^2-\pi^2m^2)^2}d\xi\\
        &=\int_0^1 \frac{(1-e^{\bi \epsilon\xi} + \bi(\epsilon\xi - \pi m))\epsilon^2\xi^2}{\sqrt{1-\xi^2}(\epsilon^2\xi^2-\pi^2m^2)^2}d\xi + (\int_{0}^{1}+\int_1^{\infty})\frac{(1-e^{-\epsilon t} - \epsilon t)\bi \epsilon^2 t^2}{\sqrt{1+t^2}(\epsilon^2 t^2 + \pi^2 m^2)^2}dt\\
        &=m^{-3}{\cal O}(\epsilon^2)+ m^{-4}{\cal O}(\epsilon^4) + \bi K^-(\epsilon;\pi m,\pi m)-\bi\int_{1}^{+\infty}\frac{\epsilon^3t^3}{\sqrt{1+t^2}(\epsilon^2t^2+\pi^2m^2)^2}dt\\
        &=\bi K^-(\epsilon;\pi m,\pi m) -\bi I_7(\epsilon) + m^{-3}{\cal O}(\epsilon^2).
\end{align*}
where 
\begin{align*}
  I_7(\epsilon) :=& \int_{1}^{+\infty}\frac{\epsilon^3t^3}{\sqrt{1+t^2}(\epsilon^2t^2+\pi^2m^2)^2}dt\\
  =&m^{-4}{\cal O}(\epsilon^3)+ \int_{2}^{+\infty}\frac{\epsilon^3t^3}{t(\epsilon^2t^2+\pi^2m^2)^2}\sum_{i=0}^{\infty}{-1/2 \choose i}t^{-2i}dt\\
  =&\sum_{i=0}^{\infty}{-1/2 \choose i}\int_{2}^{+\infty}\frac{\epsilon^3t^{2-2i}}{(\epsilon^2t^2+\pi^2m^2)^2}dt\\
  =&8\int_{\epsilon}^{+\infty}\frac{t^2}{(4t^2+\pi^2m^2)^2}dt + {\cal O}(\epsilon^2)m^{-3}.\\
  =&\frac{1}{4m} + {\cal O}(\epsilon^2)m^{-3}.
\end{align*}
Consequently, we get
\[
  d_{mm} = \bi C_0^{-}(\pi m,\pi m) -\frac{\bi}{4m} + {\cal O}(\epsilon^2)m^{-3}.
\]
Finally, when $0\neq m=n\nmid 2$, one similarly gets
\begin{align*}
  d_{mm}&=\frac{1}{2}\int_{-\infty}^{\infty}\frac{\epsilon^2\xi^2(1+\cos(\epsilon\xi))}{\sqrt{1-\xi^2}(\epsilon^2\xi^2-\pi^2m^2)^2}d\xi\\
        &=\bi K^+(\epsilon;\pi m,\pi m) -\bi I_7(\epsilon) + m^{-3}{\cal O}(\epsilon^2)\\
        &=\bi C_0^+(\pi m,\pi m) -\frac{\bi}{4m}  + m^{-3}{\cal O}(\epsilon^2\log\epsilon).
\end{align*}
\section*{Acknowledgment}
W.L. would like to thank Prof. Hai Zhang of Hong Kong University of Science and
Technology, and Prof. Junshan Lin of Arburn University for inspiring this work,
and also would like to thank Prof. Ory Schnitzer of Imperial College London for
sharing me a complete version of \cite{holsch19}.
\bibliographystyle{plain}
\bibliography{wt}

\end{document}